\documentclass[preprint,1p,number,sort&compress]{elsarticle}

\usepackage{amsmath,amsfonts,amsthm}

\usepackage[algo2e,ruled,vlined,linesnumbered]{algorithm2e}
\usepackage{hhline}

\usepackage{pgfplots}
\pgfplotsset{compat=newest}
\usetikzlibrary{plotmarks}
\usetikzlibrary{arrows.meta}
\usepgfplotslibrary{patchplots}
\usepackage{grffile}

\usepackage{listings}

\usepackage{hyperref}

\newcommand{\ee}{\mathrm{e}}
\newcommand{\RR}{\mathbb{R}}

\begin{document}
\myfooter[L]{}
\begin{frontmatter}
\title{A physics-based adaptive low-rank method for 6D Vlasov simulations}
\author[1]{Fabio Cassini\corref{cor1}}
\ead{fabio.cassini@unipd.it}
\author[2]{Lukas Einkemmer}
\ead{lukas.einkemmer@uibk.ac.at}

\cortext[cor1]{Corresponding author}
\affiliation[1]{organization={Department of Mathematics ``Tullio Levi-Civita'',
    University of Padua},addressline={Via Trieste, 63},
  postcode={35131},
  city={Padova},
country={Italy}}
\affiliation[2]{organization={Department of Mathematics,
    University of Innsbruck},addressline={Technikerstrasse, 13},
  postcode={A-6020},
  city={Innsbruck},
country={Austria}}

\begin{abstract}
The numerical simulation of kinetic plasma models in six-dimensional phase space
poses significant computational and memory challenges, particularly when
considering modern high-performance computing architectures. In this work, we propose an adaptive
dynamical low-rank approach for the efficient solution of the 6D
Vlasov--Poisson system. Our strategy is an extension of the projector-splitting
dynamical low-rank integrator. The algorithm changes the rank based on physical
observables (e.g., the electric field and energy) as opposed to controlling the
error in the distribution function. Our approach is memory efficient in the
sense that only rank $r$ objects are required to implement the scheme. This
makes the algorithm suitable for performing high-dimensional computations on
graphic processing units. Comparisons with rank adaptive Basis-Update
\& Galerkin integrators exhibit both significant memory savings and reduction in
run time. We show, by presenting a variety of 6D benchmark experiments
implemented in the C++/CUDA low-rank framework Ensign, that our technique allows us to
run simulations with a significantly smaller rank, while still preserving
the physical fidelity. Thus, the approach substantially
reduces memory footprint and computational cost and enables simulations that
are prohibitively expensive for full-rank or less memory-efficient low-rank schemes. 
\end{abstract}

\begin{keyword}
Dynamical low-rank \sep
rank adaptivity \sep
Vlasov--Poisson \sep
Plasma physics \sep
GPU-accelerated simulations
\MSC[2020]
65M22 \sep
65Y05 \sep
35Q83 \sep
65F55
\end{keyword}

\end{frontmatter}

\section{Introduction}\label{sec:intro}

Evolutionary partial differential equations (PDEs) provide a useful
framework for mathematically modeling physical phenomena. In particular, kinetic equations are widely used in plasma physics due to their ability to accurately capture fundamental phenomena~\cite{B13}. In this paper we will consider the prototypical example of a kinetic model in plasma physics, namely the Vlasov--Poisson equation
\begin{equation} \label{eq:vp}
  \left\{\begin{aligned}
    &\partial_t f(t,x,v) + v \cdot \nabla_x f(t,x,v) - E(f)(t,x) \cdot \nabla_v f(t,x,v) = 0, \\
    &E(f)(t,x) = -\nabla_x \phi(t,x), \\
    & -\Delta\phi(t,x) = \rho(f)(t,x) + 1, \quad \rho(f)(t,x) = -\int_{\Omega_v}f(t,x,v) \,dv.
  \end{aligned}
  \right.
\end{equation}
Here $f(t,x,v)$ is the particle-density function of the species under
consideration (in our case electrons), $t\in [0,T]$ represents the time,
while $x\in\Omega_x \subset \RR^d$ and $v\in\Omega_v \subset \RR^d$
refer to the spatial and velocity variables, respectively, for $d=1,2,3$.
Depending on the physical phenomenon under investigation, system~\eqref{eq:vp} is
complemented by appropriate initial and boundary conditions.

Despite the utility of kinetic equations, exact analytical solutions are rarely available due to
the inherent mathematical complexity of these models.
Developing reliable and efficient numerical methods capable of yielding accurate
solutions is therefore of paramount importance~\cite{DP14}.
In this context, a major bottleneck in simulating complex physical phenomena is
the so-called \textit{curse of dimensionality}, that is, the vast number of degrees of
freedom in the discretized system renders standard integrators
computationally prohibitive in terms of memory demands and execution times.
This is particularly critical in kinetic modeling with three-dimensional spatial
and velocity domains (i.e., $d=3$), which yield a six-dimensional phase space.
Such problems remain feasible only, if at all, on large clusters, i.e.,
on supercomputers (see, e.g, \cite{burau2010, kormann2019,einkemmer2023semi}).

A recent yet promising approach to mitigate this challenge are low-rank
approximations of the phase space (see, e.g.,~\cite{EKKMQ25} for a review).
We consider the dynamical low-rank (DLR) approach here, see \cite{KL07} or Section 2 for more details.
In this framework, the distribution function is approximated on a rank-$r$
manifold $\mathcal{M}$ (with $r$ much lower than the number of degrees of freedom)
by constraining the dynamics to the
tangent space of $\mathcal{M}$. The primary advantage of this methodology is that
computations are performed in a significantly lower-dimensional space, while
preserving the fundamental features of the system. A non-exhaustive list of contexts in which the DLR approach has been successfully applied includes matrix DEs~\cite{HNS23,CEKL24,LO14c},
the Vlasov equation~\cite{EL18,EOP20,CE22}, radiation transport~\cite{PMF20},
biochemistry~\cite{PEL23,EMP24}, the Schr\"{o}dinger equation~\cite{HNS23,CL22},
and uncertainty quantification~\cite{KCEF22}.

In the simplest case, the rank of the manifold $\mathcal{M}$  is kept
constant throughout the entire time integration. This facilitates
theoretical analysis and reduces coding effort. However, a fixed-rank approach
is clearly suboptimal in general. In fact, the minimal rank required to achieve
a prescribed accuracy is rarely known a priori and can vary significantly over
time as the system evolves (for instance, during transitions between linear
and nonlinear regimes). In this context, \textit{rank-adaptive} DLR 
algorithms offer a clear advantage over fixed-rank strategies by using an error
estimate which dynamically
adjusts the rank to meet accuracy requirements without wasting computational
resources. Such DLR methods have been derived and studied, e.g., in~\cite{ceruti2022rank,CLS23,CKL24,HW22,HNS23,hauck2023}.

For low-rank methods that construct intermediate approximations of rank higher than
$r$, rank adaptivity is natural in the sense that instead of truncating back to
a fixed rank $r$, they can truncate to a certain tolerance. This is the case for
both step-and-truncate type methods (see, e.g., \cite{ehrlacher2017,guo2022})
and the so-called rank-adaptive Basis-Update \& Galerkin (BUG) dynamical low-rank
integrators (see, e.g., \cite{ceruti2022rank,CLS23}, although not for the
parallel BUG integrator in \cite{CKL24}).
The disadvantage of such methods is that intermediate approximations of higher
rank increase the amount of memory required to run such simulations. This is
especially an issue on the more memory constraint graphic processing unit
(GPU) based systems now commonly used in high-performance computing (HPC)
systems.

To the best of our knowledge, all of the rank adaptive methods studied in the
literature use the singular values of the solution matrix to determine whether
the rank should be increased or decreased. When solving a kinetic equation, this
implies that the error in the $L^2$ norm of the distribution function is
kept below a certain tolerance. This is usually a very stringent requirement
as the distribution function can have small scale oscillations. Thus, even a
small phase error can result in errors of the order $\mathcal{O}(1)$ in this
metric. In practical simulations one is usually much more interested in averaged
quantities such as the density, the electric field, the electric energy,
or higher-order moments of the distribution function. The rank required to
effectively resolve these physical observables is commonly significantly smaller
than what would be required to control the error in the distribution function
(see, e.g., the four-dimensional numerical experiments reported in \cite{EKKMQ25}).

The present paper addresses these limitations by
\begin{itemize}
  \item Developing a rank-adaptive dynamical low-rank algorithm based on the pro\-jec\-tor-splitting (PS) integrator that never forms any objects that are larger than rank $r$. This algorithm has been implemented in the low-rank framework Ensign for both multi-core CPU and GPU based systems and we demonstrate its effectiveness for a number of six-dimensional Vlasov--Poisson simulations.
  \item We use the error in physical observables to choose the rank of the simulation. This results in a significant reduction of the rank required especially in the nonlinear regime (e.g. nonlinear Landau damping, two-stream instability, etc.). We show that while this, in general, leads to a larger error in the distribution function, the main physical features of the solution are still captured using a significantly reduced computational effort. We duly note that using such an error estimates is not restricted to our PS based approach, but could also be used in other low-rank schemes.
  \item To the best of our knowledge we provide the first direct comparison of first- and second-order PS and rank-adaptive BUG integrators for high-dimensional problem solved by an optimized code. We find that the PS integrator realizes the expected memory savings (which allows us to solve much larger problems in the GPU implementation) and reduces the run time by approximately 20\%. 
\end{itemize}

The remainder of the paper is structured as follows. In
Section~\ref{sec:recall} we review the fundamental ideas of DLR
approaches. In Section~\ref{sec:novel} we present the new procedure and
its HPC implementation in the DLR framework Ensign. In
Section~\ref{sec:numexp}, we present a collection of 6D examples
that demonstrate the efficiency and effectiveness of the approach compared
with the state of the art. Finally, we draw conclusions in Section~\ref{sec:concl}.

\section{Review of the projector-splitting DLR method for Vlasov--Poisson}\label{sec:recall}

In this section, we briefly recall how to obtain a fixed-rank PS 
DLR
approximation for the Vlasov--Poisson system~\eqref{eq:vp}. 
We invite interested readers to consult, e.g., Reference~\cite{EL18} for a
thorough explanation with full details. 
We denote the phase space by $\Omega = \Omega_x \times \Omega_v$.
Thus, rather than directly solving system~\eqref{eq:vp}, we seek an
approximation of the particle-density function $f(t,x,v)$ that, for fixed $t$,
belongs to the rank-$r$ manifold
{\small
  \begin{equation*}
  \begin{aligned}
    \mathcal{M} =
    \Big\{g(x,v) \in L^2(\Omega) : g(x,v) = \sum_{i,j}X_i(x)S_{ij}V_j(v)
    \text{ with invertible } S=(S_{ij})\in\RR^{r\times r}, \\
    X_i \in L^2(\Omega_x), V_j \in L^2(\Omega_v) \text{ with }
    \langle X_i,X_k\rangle_x =\delta_{ik}, \langle V_j,V_\ell\rangle_v =\delta_{j\ell} \Big\}.
  \end{aligned}
\end{equation*}
}
Notice that the tangent space of $\mathcal{M}$ is given by
{\small
\begin{equation*}
  \begin{aligned}
    \mathcal{T}_f\mathcal{M} =
    \Big\{\dot{g}(x,v) \in L^2(\Omega) : \dot{g}(x,v) = \sum_{i,j}\big(X_i(x)\dot{S}_{ij}V_j(v)
    + \dot{X}_i(x)S_{ij}V_j(v) + X_i(x)S_{ij}\dot{V}_j(v)\big), \\
    \text{with } \dot{S} \in \RR^{r\times r},
     \dot{X}_i \in L^2(\Omega_x), \dot{V}_j \in L^2(\Omega_v)
    \text{ and } \langle X_i,\dot{X}_k\rangle_x =0, \langle V_j,\dot{V}_\ell\rangle_v =0 \Big\},
  \end{aligned}
\end{equation*}
}%
where $\langle \cdot,\cdot\rangle_x$ and $\langle \cdot,\cdot\rangle_v$ are
the classical inner products on $L^2(\Omega_x)$ and $L^2(\Omega_v)$, respectively,
while the upper dot denotes the time derivative.
Furthermore, unless otherwise specified, all indices run from 1 to $r$ and their
explicit limits are omitted for notational simplicity.
A PS algorithm is obtained by explicitly computing
\begin{equation}\label{eq:tbproj}
  \partial_t f(t,x,v) = - P(f) (v \cdot \nabla_x f(t,x,v) -E(f) \cdot \nabla_v f(t,x,v)).
\end{equation}
Here, $P(f)$ represents the orthogonal projector onto the tangent space
$\mathcal{T}_f\mathcal{M}$, and for simplicity we keep using the symbol $f$
to denote the low-rank \textit{approximation} to the particle-density function.
Note that in the DLR framework
  the dynamics is constrained at a \emph{continuous} level to be of low rank.
  This is different from classical low-rank, where we
  leave the rank-$r$ manifold at each numerical step and then project on $\mathcal{M}$
  before starting a new step (thus resulting in larger intermediate costs, in general).
We then proceed by inserting the equivalence
\begin{equation*}
  P(f)g = P_{\overline{V}}g - P_{\overline{V}}P_{\overline{X}}g + P_{\overline{X}}g
\end{equation*}
in equation~\eqref{eq:tbproj},
where $g$ denotes a generic function and
$P_{\overline{X}}$ and $P_{\overline{V}}$ are orthogonal projectors onto
the spaces spanned by the functions $X_i$ and $V_j$, respectively
(see Reference~\cite{LO14c}).
By naturally splitting the three resulting terms, we obtain
\begin{align}
  \partial_t f(t,x,v) =
  -P_{\overline{V}}(v \cdot \nabla_x f(t,x,v)
  -E(f) \cdot \nabla_v f(t,x,v)), \label{eq:alKst}\\
  \partial_t f(t,x,v) =
  P_{\overline{V}}P_{\overline{X}}(v \cdot \nabla_x f(t,x,v)
  -E(f) \cdot \nabla_v f(t,x,v)), \label{eq:alSst}\\
  \partial_t f(t,x,v) =
  -P_{\overline{X}}(v \cdot \nabla_x f(t,x,v)
  -E(f) \cdot \nabla_v f(t,x,v)). \label{eq:alLst}
\end{align}
We are now left with explicitly applying the projectors $P_{\overline{V}}$ 
and $P_{\overline{X}}$. Starting with equation~\eqref{eq:alKst},
we define the so-called \textit{C coefficients}
\begin{equation}\label{eq:Ccoeffs}
  c_{j\ell}^1 = \int_{\Omega_v}v\tilde{V}_j(v)\tilde{V}_\ell(v)dv,
  \quad c_{j\ell}^2 = \int_{\Omega_v}\tilde{V}_j(v)\nabla_v\tilde{V}_\ell(v)dv
\end{equation}
and write the approximate particle-density function as
\begin{equation}\label{eq:Kdef}
  f(t,x,v) = \sum_{j} K_j(t,x)V_j(t,v),
  \quad K_j(t,x) = \sum_{i} X_i(t,x)S_{ij}(t).
\end{equation}
Then, straightforward calculations~\cite{EL18} lead to
\begin{equation}\label{eq:Kstep}
  \left\{
  \begin{aligned}
    &\partial_tK_j(t,x) = -\sum_{\ell}c_{j\ell}^1\cdot\nabla_xK_\ell(t,x) +
    \sum_{\ell}c_{j\ell}^2\cdot E(K)(t,x)K_\ell(t,x), \\
    &V_j(t,v) =\tilde{V}_j(v),
  \end{aligned}
  \right.
\end{equation}
which is commonly referred to as the \emph{K step} of the
PS DLR algorithm. Notice that here the velocity-dependent
factors $V_j$ remain constant in time.

Proceeding with equation~\eqref{eq:alSst} and applying both projectors
$P_{\overline{V}}$ and $P_{\overline{X}}$,
by introducing the so-called \textit{D coefficients}
\begin{equation}\label{eq:Dcoeffs}
  d_{ik}^1[E(S)] = \int_{\Omega_x}\breve{X}_i(x)E(S)\breve{X}_k(x)dx,
  \quad d_{ik}^2 = \int_{\Omega_x}\breve{X}_i(x)\nabla_x\breve{X}_k(x)dx,
\end{equation}
we get the system of ODEs for the $S$ factor
 \begin{equation}\label{eq:Sstep}
  \left\{
  \begin{aligned}
    &\dot{S}_{ij}(t) = \sum_{k,\ell}(c_{j\ell}^1\cdot d_{ik}^2-c_{j\ell}^2\cdot d_{ik}^1[E(S)])S_{k\ell}(t), \\
    &X_i(t,x)=\breve{X}_i(x), \\
    &V_j(t,v)=\tilde{V}_j(v).
  \end{aligned}
  \right.
 \end{equation}
This is usually referred to as the \emph{S step} of the
procedure. Note that at this stage neither $X_i$ nor $V_j$
change in time. 

Finally, expressing the particle-density function as
\begin{equation}\label{eq:Ldef}
  f(t,x,v) = \sum_{i} X_i(t,x)L_i(t,v), \quad L_i(t,v) = \sum_{j} S_{ij}(t)V_j(t,v),
\end{equation}
it follows that solving equation~\eqref{eq:alLst} is equivalent to
\begin{equation}\label{eq:Lstep}
  \left\{
  \begin{aligned}
    &\partial_tL_i(t,v) = \sum_{k}d_{ik}^1[E(L)]\cdot\nabla_vL_k(t,v) - \sum_{k}(d_{ik}^2\cdot v)L_k(t,v), \\
    &X_i(t,x)=\breve{X}_i(x). \\
  \end{aligned}
  \right.
\end{equation}
This stage of the PS DLR algorithm constitutes the so-called
\emph{L step}, in which the factors $X_i$ do not vary in time.

 Similarly, the electric field equations can also be formulated directly
in terms of  the low-rank factors $X_i(t,x)$, $S_{ij}(t)$ and $V_j(t,v)$,
bypassing the need to reconstruct the full distribution function $f$ (which would
obviously be impractical).
In fact, the charge density $\rho(f)(t,x)$ may take one of the following forms
\begin{align}
  &\rho(K)(t,x) = -\sum_{j} K_j(t,x)\rho(\tilde{V}_j(v)), \qquad \rho(\tilde{V}_j(v))=\int_{\Omega_v}\tilde{V}_j(v) \,dv, \label{eq:chargelrK}\\
  &\rho(S)(t,x) = -\sum_{i,j} \breve{X}_i(x)S_{ij}(t)\rho(\tilde{V}_j(v)), \nonumber \\
  &\rho(L)(t,x) = -\sum_{i} \breve{X}_i(x)\rho(L_i(t,v)), \qquad \rho(L_i(t,v)) = \int_{\Omega_v}L_i(t,v) \,dv.
  \label{eq:chargelrL}
\end{align}

An approximate solution to the Vlasov--Poisson system~\eqref{eq:vp} is
ultimately obtained by combining the partial solutions of the K, S, and L steps
within an operator-splitting framework, as detailed, e.g.,
in~\cite[Algorithms 1 and 2]{CE22}. In the simplest setting, a first-order
Lie splitting scheme concatenates the three subflows in sequence.
For second-order temporal accuracy, we may employ an almost-symmetric Strang
splitting scheme. In this case, the electric field is first evaluated at the half-step
exploiting the low-rank equations, after which the time integration proceeds via
the standard Strang splitting (see~\cite{CE22} for complete algorithmic details).
In both cases, after a phase space semidiscretization with $n=n_{x_i}=n_{v_i}$
degrees of freedom in each direction, the resulting memory and
computational complexity of the PS DLR algorithms theoretically scale as
$\mathcal{O}(rn^d)$ and $\mathcal{O}(r^2n^d)$, respectively, compared to
$\mathcal{O}(n^{2d})$ of a full-rank simulation.
We refer the reader to~\cite{EL18} for a more detailed complexity analysis.

\begin{algorithm2e}
  \SetAlgoLined
  \KwIn{$X^0$, $S^0$, $V^0$ s.t.~$f(0,x,v) \approx \sum_{i,j}X^0_i(x)S^0_{ij}V^0_j(v)$, time step size $\tau$}
  \KwOut{$X^\tau$, $S^\tau$, $V^\tau$ s.t. $f(\tau,x,v) \approx \sum_{i,j}X^\tau_i(x)S^\tau_{ij}V^\tau_j(v)$}
  Compute C coefficients~\eqref{eq:Ccoeffs} using $V^0$\;
  Compute $K^0$ using $X^0$ and $S^0$ as in equation~\eqref{eq:Kdef}\;
  Compute the electric field $E^0$ using $K^0$, $V^0$, and equation~\eqref{eq:chargelrK}\;
  Solve equation~\eqref{eq:Kstep} with initial value $K^0$ and $E^0$ up to $\tau$ to obtain $K^1$\;
  Compute a QR decomposition of $K^1$ to obtain $X^1$ and $S^1$\;
  Compute D coefficients~\eqref{eq:Dcoeffs} using $E^0$ and $X^1$\;
  Solve equation~\eqref{eq:Sstep} with initial value $S^1$ and $E^0$ up to $\tau$ to obtain $S^2$\;
  Compute $L^0$ using $V^0$ and $S^2$ as in equation~\eqref{eq:Ldef}\;
  Solve equation~\eqref{eq:Lstep} with initial value $L^0$ and $E^0$ up to $\tau$ to obtain $L^1$\;
  Compute a QR decomposition of $L^1$ to obtain $V^1$ and $S^3$\;
  $X^\tau \gets X^1$,  $S^\tau \gets S^3$, and  $V^\tau \gets V^1$\;
  \caption{Excerpt of a step of a first-order PS DLR integrator for the Vlasov--Poisson equations~\eqref{eq:vp}.}
  \label{alg:firstorder}
\end{algorithm2e}
For the reader's convenience and for later reference, we provide the pseudocode
for the first- and second-order methods in Algorithms~\ref{alg:firstorder} and~\ref{alg:secondorder},
assuming a phase-space discretization that collects the degrees of freedom in suitable matrices.
\begin{algorithm2e}
  \SetAlgoLined
  \KwIn{$X^0$, $S^0$, $V^0$ s.t.~$f(0,x,v) \approx \sum_{i,j}X^0_i(x)S^0_{ij}V^0_j(v)$, time step size $\tau$}
  \KwOut{$X^\tau$, $S^\tau$, $V^\tau$ s.t. $f(\tau,x,v) \approx \sum_{i,j}X^\tau_i(x)S^\tau_{ij}V^\tau_j(v)$}
  Perform steps 1--9 of Algorithm~\ref{alg:firstorder} with time step size $\tau/2$\;
  Compute the electric field $E^{1/2}$ using $X^1$ and $L^1$ from step 1 and equation~\eqref{eq:chargelrL}\;
  Solve equation~\eqref{eq:Kstep} with initial value $K^0$ and $E^{1/2}$ up to $\tfrac{\tau}{2}$ to obtain $K^2$\;
  Compute a QR decomposition of $K^2$ to obtain $X^2$ and $S^3$\;
  Compute D coefficients~\eqref{eq:Dcoeffs} using $E^{1/2}$ and $X^2$\;
  Solve equation~\eqref{eq:Sstep} with initial value $S^3$ and $E^{1/2}$ up to $\tfrac{\tau}{2}$ to obtain $S^4$\;
  Compute $L^1$ using $V^0$ and $S^4$ as in equation~\eqref{eq:Ldef}\;
  Solve equation~\eqref{eq:Lstep} with initial value $L^1$ and $E^{1/2}$ up to $\tau$ to obtain $L^2$\;
  Compute a QR decomposition of $L^2$ to obtain $V^1$ and $S^5$\;
  Compute C coefficients~\eqref{eq:Ccoeffs} using $V^1$\;
  Solve equation~\eqref{eq:Sstep} with initial value $S^5$ and $E^{1/2}$ up to $\tfrac{\tau}{2}$ to obtain $S^6$\;
  Compute $K^3$ using $X^2$ and $S^6$ as in equation~\eqref{eq:Kdef}\;
  Solve equation~\eqref{eq:Kstep} with initial value $K^3$ and $E^{1/2}$ up to $\tfrac{\tau}{2}$ to obtain $K^4$\;
  Compute a QR decomposition of $K^4$ to obtain $X^3$ and $S^7$\;
  $X^\tau \gets X^3$,  $S^\tau \gets S^7$, and  $V^\tau \gets V^1$\;
  \caption{Excerpt of a step of a second-order PS DLR integrator for the Vlasov--Poisson equations~\eqref{eq:vp}.}
  \label{alg:secondorder}
\end{algorithm2e} 

\section{The novel rank-adaptive procedure and its Ensign implementation}\label{sec:novel}
As previously noted, in DLR algorithms selecting an appropriate rank $r$ is of
great importance for practical applications, particularly in large-scale simulations.
In fact, this choice involves a fundamental trade-off. On one hand, a higher rank
is desirable to maintain an appropriate accuracy and correctly capture the physical
dynamics of the underlying model. On the other hand, keeping an excessively
high rank throughout the entire simulation leads to prohibitive computational
costs (even when far below the full-rank limit in a 6D phase space setting).
Although fixed-rank strategies have been extensively studied and applied in
many scenarios, a rank-adaptive approach represents a far more natural alternative.
Indeed, as the system evolves in time, transient multiscale phenomena would 
naturally modify the required approximation subspace, requiring different
ranks in different time intervals. For instance, in a Vlasov simulation of a
classical two-stream instability, the initial exponential growth of electric energy
evolves to a nonlinear saturation phase, with each regime necessitating a
different rank to remain both accurate and efficient
(see the next section for the numerical simulations).

As mentioned in the introduction, while rank adaptation is not new in the
DLR literature, to the best of our knowledge this work introduces
the first HPC-oriented strategies that base rank selection directly on physical
observables of the PDE rather than purely matrix quantities.
The idea is that guiding the rank choice through the specific physics of the
system provides a more targeted approximation rather than treating the model as
a generic matrix differential equation. Because our focus is on
Vlasov simulations, we consider three physical quantities:
\begin{itemize}
  \item The distribution function (approach labeled as \texttt{df}).
  \item The electric field (approach labeled as \texttt{ef}).
  \item The electric energy (approach labeled as \texttt{ee}).
\end{itemize}
Each of these three strategies is described in detail below.
A single pseudocode summarizing the procedures is presented in Algorithm~\ref{alg:radapt}.

\begin{algorithm2e}[!htb]
  \SetAlgoLined
  \KwIn{rank $r$, $X^0$, $S^0$, $V^0$ s.t~$f(0,x,v) \approx \sum_{i,j}^rX^0_i(x)S^0_{ij}V^0_j(v)$,
  time step size $\tau$, tolerances $\epsilon_{\mathrm{decr}}$ and $\epsilon_{\mathrm{incr}}$ (absolute and relative),
  adaptivity control $\kappa$ (\texttt{df}, \texttt{ef}, or \texttt{ee})}
  \KwOut{rank $r_{\mathrm{new}}$, $X^\tau$, $S^\tau$, $V^\tau$ s.t.~$f(0,x,v) \approx \sum_{i,j}^{r_{\mathrm{new}}}X^\tau_i(x)S^\tau_{ij}V^\tau_j(v)$}
  Perform steps of Algorithm~\ref{alg:firstorder} or~\ref{alg:secondorder}
  to get the outputs $X^\tau$, $S^\tau$, and $V^\tau$\;
  Compute SVD $S^\tau = A\Lambda B^{\sf T}$, $\Lambda$ with singular values $\lambda_k$ in descending order\;
  \uIf{$\kappa == \texttt{ef}$ or $\kappa == \texttt{ee}$}
  {
    Compute the electric field $E$\;
    Form $\Lambda_{\mathrm{red}}$ setting $\lambda_r$ to $0$\;
    $S^\tau_{\mathrm{red}} \gets A\Lambda_{\mathrm{red}} B^{\sf T}$\;
    Compute the reduced electric field $E_{\mathrm{red}}$\;
    \uIf{$\kappa == \texttt{ef}$}
    {
      $\alpha \gets \lVert E-E_{\mathrm{red}} \rVert_{L^2}$ and $\ell \gets \lVert E \rVert_{L^2}$ // \texttt{ef approach}\;
    }
    \Else
    {
      Compute the electric energy $\delta$ and reduced electric energy $\delta_{\mathrm{red}}$\;
      $\alpha \gets |\delta - \delta_{\mathrm{red}}|$ and $\ell \gets \delta$ // \texttt{ee approach}\;
    }
  }
  \Else
  {
    $\alpha \gets \lambda_r$ and $\ell \gets 0$ // \texttt{df approach}\;
  }
  \uIf{$\alpha\leq\epsilon_{\mathrm{decr}}^{\mathrm{abs}}+\ell\epsilon_{\mathrm{decr}}^{\mathrm{rel}}$}
  {
    $r_{\mathrm{new}} \gets r-1$ \texttt{// Accepted step, decreasing rank}\;
    $X^\tau \gets X^\tau(:,1:r_{\mathrm{new}})$, $S^\tau \gets S^\tau(1:r_{\mathrm{new}},1:r_{\mathrm{new}})$, $V^\tau \gets V^\tau(:,1:r_{\mathrm{new}})$\;
  }
  \uElseIf{$\alpha\geq\epsilon_{\mathrm{incr}}^{\mathrm{abs}} + \ell\epsilon_{\mathrm{incr}}^{\mathrm{rel}}$}
  {
    $r_{\mathrm{new}} \gets r+1$ \texttt{// Rejected step, increasing rank}\;
    Add a random column to $X^0$ and to $V^0$, a null row and column to $S^0$\;
    Compute a QR decomposition of $X^0$ to obtain $X^0_{\mathrm{new}}$ and $G$\;
    Compute a QR decomposition of $V^0$ to obtain $V^0_{\mathrm{new}}$ and $H$\;
    $S^0_{\mathrm{new}} \gets GS^0H^{\sf T}$\;
    $X^0 \gets X^0_{\mathrm{new}}$, $S^0 \gets S^0_{\mathrm{new}}$, and $V^0 \gets S^0_{\mathrm{new}}$\;
    Redo the procedure going to step 1\;
  }
  \Else
  {
    $r_{\mathrm{new}} \gets r$ \texttt{// Accepted step, maintaining rank}\;
  }
  \caption{The proposed rank-adaptive PS algorithm.}
  \label{alg:radapt}
\end{algorithm2e}

\subsection{Distribution function approach \texttt{df}}

Selecting adaptively the rank by using information directly from the
distribution function is the closest approach to established methods in the
literature for generic matrix differential equations~\cite{HNS23,ceruti2022rank}.
Specifically, this technique relies on tracking the singular values of the
matrix representing the discretized distribution function at each computational
node. In this case, the magnitude of the available singular values
measures the amount of information present at a given time step.

More specifically, assume that at a fixed time the numerical solution $f(t,x,v)$
on the computational grid is represented by a matrix $F\in\RR^{N_x \times N_v}$,
where $N_x=\prod n_{x_i}$ and $N_v=\prod n_{v_i}$ denote the spatial and velocity degrees of freedom,
respectively. Notice that in the DLR framework $F$ is available via its
low-rank factorization $F = X S V^{\sf T}$, where $X \in \RR^{N_x \times r}$,
$S\in \RR^{r \times r}$, and $V \in \RR^{N_v \times r}$ (see Section~\ref{sec:recall}).
Consequently, the singular values of $F$ can be easily computed from the
small-sized coefficient matrix $S$ (e.g., through a classical singular value
decomposition).
Let $\lambda_1 \geq \lambda_2 \geq \ldots \geq \lambda_r$ denote the
singular values (in descending order). The smallest 
singular value available, that is $\lambda_r$, can then be used to decide
whether to increase, decrease, or maintain the current rank, comparing its
magnitude against some pre-defined user tolerances. In fact, we perform one among the following choices:
\begin{itemize}
  \item Rank \textit{reduction}: if $\lambda_r\leq\epsilon_{\mathrm{decr}}$,
  the smallest available singular value relates to negligible information. Therefore,
  the current step is accepted, and the rank is decremented for the subsequent step
  by truncating the corresponding rows and columns of the low-rank factors
  $X$, $S$, and $V$.
  \item Rank \textit{increment}: if $\lambda_r\geq\epsilon_{\mathrm{incr}}$,
  the smallest available singular value is still significant. Therefore,
  the current step is rejected, the rank is incremented, and the time stepping
  procedure is restarted. In practice, this is achieved by expanding the relevant
  matrices with random columns followed by an orthonormalization step.
  \item Rank \textit{maintenance}: if $\epsilon_{\mathrm{decr}} < \lambda_r < \epsilon_{\mathrm{incr}}$,
  the current rank $r$ is adequate; the step is accepted, and the time integration
  proceeds maintaining the same rank.
\end{itemize}

Because the just described distribution function-based approach proceeds independently
of the underlying PDEs, we will employ it as a natural benchmark.
This criterion is closely related to the truncation-error criterion used by
rank-adaptive BUG integrators~\cite{ceruti2022rank,CEKL24}. Note, however,
that since the BUG integrator increases the intermediate
rank to $2r$, it only performs rank reduction and therefore uses a single tolerance.

In kinetic Vlasov simulations,
however, physical observables such as the electric field and electric energy are often of primary interest. 
We thus propose two alternative rank-adaptive strategies based directly on these quantities.
This generally results in a lower rank and thus reduced computational burden (see the numerical experiments in Section~\ref{sec:numexp}).

\subsection{Electric field approach \texttt{ef}}
At the end of each time step, we may choose to adapt the rank based on information
coming from the electric field. To this aim, we first compute the electric field
$E$. To estimate the truncation error, we compute a \textit{reduced} electric field
$E_{\mathrm{red}}$ by omitting the information carried  by the smallest available singular value $\lambda_r$.
To this aim, similarly to the previous approach, we compute a singular value decomposition
of the coefficient matrix $S$ and remove the corresponding rows and columns in the
relevant matrices.
The discrete $L^2$ norm difference between the full and reduced electric fields
yields an error estimate $\alpha = \lVert E- E_{\mathrm{red}} \rVert_{L^2}$.
We compare it with separate mixed absolute-relative tolerances for rank decrease
and increase, using $\ell=\lVert E\rVert_{L^2}$ as the scale:
\begin{itemize}
  \item Rank \textit{reduction}: if
  $\alpha\leq\epsilon_{\mathrm{decr}}^{\mathrm{abs}}+
  \ell\epsilon_{\mathrm{decr}}^{\mathrm{rel}}$, the current
  step is accepted, and the rank is decremented for the subsequent time step.
  \item Rank \textit{increment}: if
  $\alpha\geq\epsilon_{\mathrm{incr}}^{\mathrm{abs}}+
  \ell\epsilon_{\mathrm{incr}}^{\mathrm{rel}}$, the current rank is insufficient to
  accurately capture the electric field. The step is rejected, the rank is
  increased, and the time integration step is restarted.
  \item Rank \textit{maintenance}: if the error lies between the two thresholds,
  the step is accepted and the rank is held constant for the following step.
\end{itemize}

\subsection{Electric energy approach \texttt{ee}}
Rank adaptation via electric energy arguments follows a similar logic.
Rather than evaluating the $L^2$ norm error of the electric fields themselves,
we measure the error directly in the electric energy. Specifically, the truncation error
is estimated by comparing the full electric energy to the reduced energy computed
from the electric field with its smallest singular value omitted.
Step acceptance and rank updates (reduction, increment, or maintenance of 
$r$) then proceed according to analogous mixed absolute-relative threshold criteria
used for the \texttt{ef} approach.

\subsection{Practical considerations}
As highlighted in Algorithm~\ref{alg:radapt},
in our implementation we decide to restrict rank updates
to adjustments by at most one per time step. On one hand,
this strategy serves as a practical safeguard against strong rank variations
between consecutive time steps. On the other hand, empirical evidence across all
our numerical experiments shows that this strategy is both robust and effective.
Naturally, the idea can be straightforwardly generalized to allow arbitrary
rank variations without extra conceptual complexity, provided it remains
computationally feasible. For instance, alternative frameworks such as the parallel
BUG integrator~\cite{CKL24} exploit weighted singular value averages
to determine rank truncations or expansions.
Additionally, we incorporate in our code explicit safety bounds on the manifold
rank $r$. Specifically, a minimum rank $r_\mathrm{min}$ is imposed to guarantee
that essential physical information is always retained, while a maximum rank
$r_\mathrm{max}$ is prescribed to prevent exceeding hardware-specific memory
constraints.

Designed specifically for HPC architectures, the proposed PS DLR
solver is
built to achieve favorable parallel scalability, particularly on GPU devices.
To this aim, the practical implementation of the algorithm builds upon the
matrix-oriented Ensign framework\footnote{Freely available under MIT licence at \url{https://github.com/leinkemmer/Ensign}.}
introduced in Reference~\cite{CE22}, which we
enriched with novel routines to support rank-adaptive dynamics efficiently.
In this sense, a primary computational bottleneck in GPU-based rank-adaptive
algorithms is the overhead associated with frequent memory reallocations.
To overcome this, at the initialization stage
we pre-allocate memory buffers for large-sized objects in terms of
degrees of freedom and maximum anticipated rank (e.g.,
we reserve memory for matrices of size $N_x \times r_\mathrm{max}$
and $N_v \times r_\mathrm{max}$ for the space and velocity terms $X$ and $V$,
respectively).
Then, throughout time integration,
dynamic rank variations are managed exclusively via pointers and sub-matrices
without the need to modify the underlying memory layout.
This strategy avoids costly memory reallocations and data transfers,
ensuring optimal hardware usage in a parallel context.

\begin{figure}[!htb]
\begin{lstlisting}[basicstyle=\scriptsize\ttfamily]
/* ... */
stloc sl = stloc::device; //works on the GPU
//stloc sl = stloc::host; //works on the CPU
int rmax = 120;
int Nx = 64*64*64;
int Nv = 128*128*128;
int r = 10;
mat X({Nx,rmax},{Nx,r},sl); //allocates memory for Nx x rmax, actual size Nx x r
mat S({r,r}, sl); //matrix of size r x r 
mat V({Nv,rmax},{Nv,r},sl); //allocates memory for Nv x rmax, actual size Nv x r
/* Compute K and L, e.g., blas.matmul(X,S,K) and blas.matmul_transb(S,V,L) */
/* Do K, S, and L steps */
rnew = 11; //obtained from adaptivity procedure
X.update_shape({Nx,rnew}); //actual size Nx x rnew, NO memory re-allocated
S.resize({rnew,rnew}); //resize matrix to rnew x rnew, memory re-allocated
V.update_shape({Nv,rnew}); //actual size Nv x rnew, NO memory re-allocated
/* ... */
\end{lstlisting}
\caption{Sketch of C++/CUDA implementation for efficient memory management using the
Ensign framework in the context of the proposed adaptive approach.
The hardware is automatically selected based on the value in \texttt{sl}.
The syntax \texttt{/*...*/} indicates code not reported for simplicity
of exposition. \texttt{mat} is a shorthand notation for \texttt{multi\_array<double,2>},
the relevant Ensign structure for matrices.}\label{code:pseudo_alloc}
\end{figure}
A representative example of memory allocation and management for the quantities
involved in Ensign is illustrated in Figure~\ref{code:pseudo_alloc}.
Note that in this framework we adopt a column-major memory layout, which is
particularly well suited to the considered low-rank factors. Indeed,
these representations naturally organize the large-sized matrices $X$ and $V$
by columns, with rank adaptation primarily involving the addition or removal of
columns (representing in fact the basis).
Since $X$ and $V$ are essentially the largest objects in our representation,
we pre-allocate their memory taking into account the maximum required rank
(lines 8 and 10). Whenever the rank changes, we update the actual matrix
dimensions by adjusting the corresponding pointers,
without reallocating or moving the underlying data (lines 14 and 16).
We stress once again that the column-major layout is essential to this approach, as
it allows changes in the actual rank to be handled by modifying the number of
columns while preserving the existing data in memory.
In contrast, the matrix $S$ is much smaller than $X$ and $V$, so pre-allocation
is unnecessary (line 9). In fact, we resize $S$ as needed (line 15).
Although handling $S$ may involve more extensive memory management (since we 
are changing both rows and columns, hence data can not be preserved), its
small size makes the associated overhead negligible. This is confirmed by the
numerical experiments presented in the following section.
For brevity of exposition, and since providing the complete executable code
would require considerably more space, we refer the interested reader to the
corresponding code in the Ensign GitHub repository,
which provides a complete example
in the context of 6D Vlasov--Poisson simulation.

\section{Numerical experiments}\label{sec:numexp}
In the following we present several numerical experiments that highlight the
superiority of the proposed PS DLR approach to perform 6D Vlasov simulations on
HPC architectures. In particular, throughout the section we will consider
the following problems:
\begin{itemize}
  \item \textit{Linear Landau}. The domain is
  $\Omega = (0,4\pi)^3\times(-6,6)^3$ and the initial particle-density function is
  {\small
  \begin{equation}\label{eq:ll}
  f_0(x,v) = \frac{1}{\sqrt{(2\pi)^3}}
  \ee^{-(v_1^2+v_2^2+v_3^2)/2}
  \times(1+\alpha_1\cos(\kappa_1 x_1)+\alpha_2\cos(\kappa_2 x_2)+\alpha_3\cos(\kappa_3 x_3)).
  \end{equation}
  }
  The parameters are set to $\alpha_1 = \alpha_2 = \alpha_3 = 10^{-2}$ and 
  $\kappa_1 = \kappa_2 = \kappa_3 = \frac{1}{2}$.

  \item \textit{Two-stream instability}. The name refers to the two counter-propagating
  streams in the $v_1$ direction. The product form of the initial condition below
  yields eight Gaussian components in three-dimensional velocity space. The domain is
  $\Omega = (0,10\pi)^3\times(-9,9)^3$ and the initial particle-density function is
  {\small
  \begin{equation}\label{eq:ts}
  \begin{aligned}
  f_0(x,v) = &\frac{1}{\sqrt{(8\pi)^3}}\left(\ee^{-(v_1-\overline{v}_1)^2/2}+\ee^{-(v_1-\tilde{v}_1)^2/2}\right)
  \times\left(\ee^{-(v_2-\overline{v}_2)^2/2}+\ee^{-(v_2-\tilde{v}_2)^2/2}\right) \\
  &\times\left(\ee^{-(v_3-\overline{v}_3)^2/2}+\ee^{-(v_3-\tilde{v}_3)^2/2}\right) \\
  &\times(1+\alpha_1\cos(\kappa_1 x_1)+\alpha_2\cos(\kappa_2 x_2)+\alpha_3\cos(\kappa_3 x_3)).
  \end{aligned}
  \end{equation}
  }
  The parameters are set to
  $\alpha_1 = \alpha_2 = \alpha_3 = 10^{-3}$, $\kappa_1 = \kappa_2 = \kappa_3 = \frac{1}{5}$,
  $\overline{v}_1 = \frac{5}{2}$, $\overline{v}_2 = \overline{v}_3 = 0$,
  $\tilde{v}_1 = -\frac{5}{2}$, $\tilde{v}_2 = -\frac{9}{4}$ and $\tilde{v}_3 = -2$.

  \item \textit{Nonlinear Landau}. The domain and the initial particle-density
  function
  are the same as for the linear Landau problem~\eqref{eq:ll}, but the parameters are 
  set to $\alpha_1 = \alpha_2 = \alpha_3 = 0.3$ and
  $\kappa_1 = \kappa_2 = \kappa_3 = \frac{1}{2}$.

  \item \textit{Bump-on-tail}. The domain is $\Omega = (0,20\pi/3)^3\times(-9,9)^3$
  and the initial particle-density function is
  {\small
  \begin{equation}\label{eq:bot}
  \begin{aligned}
  f_0(x_1,x_2,x_3,v_1,v_2,v_3) = &\frac{1}{\sqrt{(2\pi)^3}}
  \ee^{-(v_2^2+v_3^2)/2}(\beta\ee^{-v_1^2/2}+\gamma\ee^{-2(v_1-\delta)^2})\\
  &\times(1+\alpha(\cos(\kappa x_1)+\cos(\kappa x_2)+\cos(\kappa x_3))).
  \end{aligned}
  \end{equation}
  }
  The parameters are set to $\alpha = \frac{3}{100}$, $\beta = \frac{9}{10}$,
  $\gamma=\frac{1}{5}$, $\delta=\frac{9}{2}$, $\kappa = \frac{3}{10}$.
\end{itemize}

All our examples use $d=3$, that is three dimensions in both space and velocity,
and periodic boundary conditions.
We use a Fourier pseudospectral approach for the discretization
of the phase space. Concerning the time evolution of the low-rank factors
(namely, in the K, S, and L steps),
we always employ a single step of the classical explicit Runge--Kutta integrator
of order four.
This choice combines satisfactory stability properties along the imaginary axis
with the absence of an impractically restrictive CFL condition in our experiments.

As mentioned earlier in the paper, for the software we employ the
C++/CUDA Ensign framework~\cite{CE22} and double precision arithmetic.
Concerning the parallelization, we use OpenMP for the CPU part
and CUDA kernels for the GPU code. The needed BLAS operations are performed
using the high-performance OpenBLAS library on the CPU~\cite{ZWZ12,QZZY13} and the NVIDIA CUDA Toolkit (specifically, cuBLAS~\cite{cublas}) on the GPU.
The Fourier pseudospectral approximations are
performed using classical FFTW-based routines~\cite{FJ05} (note that cuBLAS provides an FFTW-compatible interface).

In terms of hardware, we employ a workstation equipped with a single-socket
AMD EPYC 9554P 64-Core CPU (192 GB of RAM) and an NVIDIA Tesla A100
graphics card with 80 GB of dedicated memory (theoretical peak memory bandwidth
of 1555 GB/s and peak double-precision performance of 9.7 TFlops).

For comparison, we consider the state-of-the-art adaptive BUG
integrator in both its first-order formulation~\cite{CL22} and one of its second-order variants
(specifically, the BUG midpoint $3r$ variant described in~\cite{CEKL24}).
We stress once again that in our context a relevant distinction between these approaches
and the PS-type proposed here is related to their memory requirements.
Indeed, in its first-order formulation, BUG requires storage of low-rank
objects with column size $2r$, while the second-order BUG midpoint scheme
requires objects of column size $3r$. By contrast, the PS integrators require objects of column size $r$ (since
no basis augmentation is required).
While these additional memory requirements
are generally not prohibitive (especially for 1+1D or 2+2D simulations),
they can become a significant limitation in the context of six-dimensional Vlasov
simulations on GPU-based HPC systems, where the available device memory is
inherently limited.
In this sense, reducing the size of the objects that must be simultaneously
stored can be crucial for making large-scale simulations feasible.

Note also that the memory efficiency of our approach is particularly relevant when
considering the actual scale of the simulations. With our approach and implementation,
we demonstrate that six-dimensional Vlasov simulations can in fact be performed on a
single suitably
equipped workstation, without requiring a large distributed-memory HPC system.
For comparison, a standard full-rank discretization of the same problem
would require storing the full six-dimensional phase space distribution,
resulting in a number of degrees of freedom and associated memory requirements
that make such simulations impractical.

We organize the numerical experiments as follows.
We first consider a set of tests based on the linear Landau damping and
two-stream instability problems. These experiments are designed to numerically
verify the order of convergence of the proposed approach, compare its performance
with the BUG integrators both on the CPU and GPU, and investigate the effect of varying the prescribed
tolerances and rank adaptivity techniques.
We then consider the whole collection of benchmark problems, namely the ones
introduced at the beginning of this section.
These simulations employ the second-order adaptive PS integrator
on the GPU and use discretizations with a large number of degrees of freedom,
pushing the available hardware resources close to their limits.
In this set of experiments, we fix the tolerance for each simulation and
use only the electric energy to control adaptivity.
For each benchmark, we report the evolution of the electric energy together
with the corresponding rank evolution, demonstrating the ability of the
adaptive strategy to adjust the rank (and hence the computational complexity)
to correctly capture the underlying dynamics.

\subsection{Convergence rate verification}\label{sec:orders}

\begin{figure}[htb]
  \centering
%
%
\begin{tikzpicture}

\begin{axis}[%
width=1.65in,
height=1.25in,
at={(1.236in,0.823in)},
scale only axis,
xmode=log,
xmin=10,
xmax=100,
xminorticks=true,
ymode=log,
ymin=1e-10,
ymax=1e-6,
xtick = {10,20,50,100},
xticklabels = {{10},{20},{50},{100}},
xlabel = {Number of steps},
ylabel = {Error},
yminorticks=true,
axis background/.style={fill=white},
title={CPU},
legend style={legend cell align=left, align=left, at={(0.97,0.67)}, font=\tiny}
]
\addplot [color=blue, line width=1.2pt, only marks, mark=x, mark size=3pt, mark options={solid, blue}]
  table[row sep=crcr]{%
10	5.82137358203827e-07\\
20	2.91595077221111e-07\\
30	1.94512540907252e-07\\
40	1.45927306279214e-07\\
50	1.16762116635017e-07\\
60	9.73216190404409e-08\\
70	8.34881433946353e-08\\
80	7.30974605819137e-08\\
90	6.50066199679135e-08\\
100	5.85300828559457e-08\\
};
\addlegendentry{PS1}

\addplot [color=red,line width=1.2pt, only marks, mark=+,mark size=3pt, mark options={solid, red}]
  table[row sep=crcr]{%
10	5.77191121822145e-07\\
20	2.90429122916965e-07\\
30	1.93992306720087e-07\\
40	1.45681351815375e-07\\
50	1.16660592686599e-07\\
60	9.72870865593793e-08\\
70	8.34662616454696e-08\\
80	7.30839475605471e-08\\
90	6.49988321475648e-08\\
100	5.85243933581894e-08\\
};
\addlegendentry{BUG1}

\addplot [color=cyan, only marks,line width=1.2pt, mark=o,mark size=3pt, mark options={solid, cyan}]
  table[row sep=crcr]{%
10	1.40683707989079e-08\\
20	3.52371371838935e-09\\
30	1.56798172973023e-09\\
40	8.82595872709206e-10\\
50	5.65022894306861e-10\\
60	3.92383141736186e-10\\
70	2.8823709910393e-10\\
80	2.20630771715611e-10\\
90	1.74288362609949e-10\\
100	1.41151491803814e-10\\
};
\addlegendentry{PS2}

\addplot [color=magenta, only marks,line width=1.2pt, mark=star,mark size=3pt, mark options={solid, magenta}]
  table[row sep=crcr]{%
10	9.88036071039236e-09\\
20	2.475319856603e-09\\
30	1.10320054337773e-09\\
40	6.23432103956451e-10\\
50	4.0123357178645e-10\\
60	2.79341522730347e-10\\
70	2.06172489257652e-10\\
80	1.60078358404782e-10\\
90	1.286331749758e-10\\
100	1.05960277608542e-10\\
};
\addlegendentry{BUG2}

\addplot [color=black, dashed, forget plot]
  table[row sep=crcr]{%
10	5.85243933581894e-07\\
100	5.85243933581894e-08\\
};
\addplot [color=black, forget plot]
  table[row sep=crcr]{%
10	1.41151491803814e-08\\
100	1.41151491803814e-10\\
};
\addplot [color=black, forget plot]
  table[row sep=crcr]{%
10	1.05960277608542e-08\\
100	1.05960277608542e-10\\
};
\end{axis}

\end{tikzpicture}%
%
%
\begin{tikzpicture}

\begin{axis}[%
width=1.65in,
height=1.25in,
at={(1.236in,0.823in)},
scale only axis,
xmode=log,
xmin=10,
xmax=100,
xminorticks=true,
ymode=log,
ymin=1e-10,
ymax=1e-6,
xtick = {10,20,50,100},
xticklabels = {{10},{20},{50},{100}},
xlabel = {Number of steps},
ylabel = {Error},
yminorticks=true,
axis background/.style={fill=white},
title={GPU},
legend style={legend cell align=left, align=left, at={(0.97,0.67)}, font=\tiny}
]
\addplot [color=blue, line width=1.2pt, only marks, mark=triangle, mark size=3pt, mark options={solid, blue}]
  table[row sep=crcr]{%
10	5.82137358522132e-07\\
20	2.91595077061958e-07\\
30	1.94512540907252e-07\\
40	1.45927305536502e-07\\
50	1.16762114672133e-07\\
60	9.73216197831537e-08\\
70	8.34881458349771e-08\\
80	7.30974625978481e-08\\
90	6.50066231509679e-08\\
100	5.8530084977982e-08\\
};
\addlegendentry{PS1}

\addplot [color=red,line width=1.2pt, only marks, mark=asterisk,mark size=3pt, mark options={solid, red}]
  table[row sep=crcr]{%
10	5.77191110363149e-07\\
20	2.90429079362169e-07\\
30	1.93992308258563e-07\\
40	1.45681393407288e-07\\
50	1.16660646321067e-07\\
60	9.72871179124664e-08\\
70	8.34662853592259e-08\\
80	7.30839972692487e-08\\
90	6.49989223341095e-08\\
100	5.85244870991452e-08\\
};
\addlegendentry{BUG1}

\addplot [color=cyan, only marks,line width=1.2pt, mark=square,mark size=3pt, mark options={solid, cyan}]
  table[row sep=crcr]{%
10	1.40683709580606e-08\\
20	3.52371265737118e-09\\
30	1.56798130532297e-09\\
40	8.82596031861933e-10\\
50	5.65022416848682e-10\\
60	3.92382080718011e-10\\
70	2.88234499609403e-10\\
80	2.2062875578108e-10\\
90	1.74284861249974e-10\\
100	1.41149051462013e-10\\
};
\addlegendentry{PS2}

\addplot [color=magenta, only marks,line width=1.2pt, mark=diamond,mark size=3pt, mark options={solid, magenta}]
  table[row sep=crcr]{%
10	9.81287486161775e-09\\
20	2.5016082264502e-09\\
30	1.10122020600625e-09\\
40	6.2781474562734e-10\\
50	4.01101103667385e-10\\
60	2.80211557633298e-10\\
70	2.07430007997868e-10\\
80	1.60986589962009e-10\\
90	1.32422070876332e-10\\
100	1.10080211179465e-10\\
};
\addlegendentry{BUG2}

\addplot [color=black, dashed, forget plot]
  table[row sep=crcr]{%
10	5.8530084977982e-07\\
100	5.8530084977982e-08\\
};
\addplot [color=black, forget plot]
  table[row sep=crcr]{%
10	1.41149051462013e-08\\
100	1.41149051462013e-10\\
};
\addplot [color=black, forget plot]
  table[row sep=crcr]{%
10	1.10080211179465e-08\\
100	1.10080211179465e-10\\
};
\end{axis}

\end{tikzpicture}%
  \caption{Convergence rate verification for the employed integrators on CPU (left)
  and GPU (right), see Section~\ref{sec:orders}. The test is the linear Landau damping problem
  with $16^3\times 16^3$ degrees of freedom up to the final time $T=0.1$. The rank is
  fixed to 10 for all methods.}
  \label{fig:orders}
\end{figure}
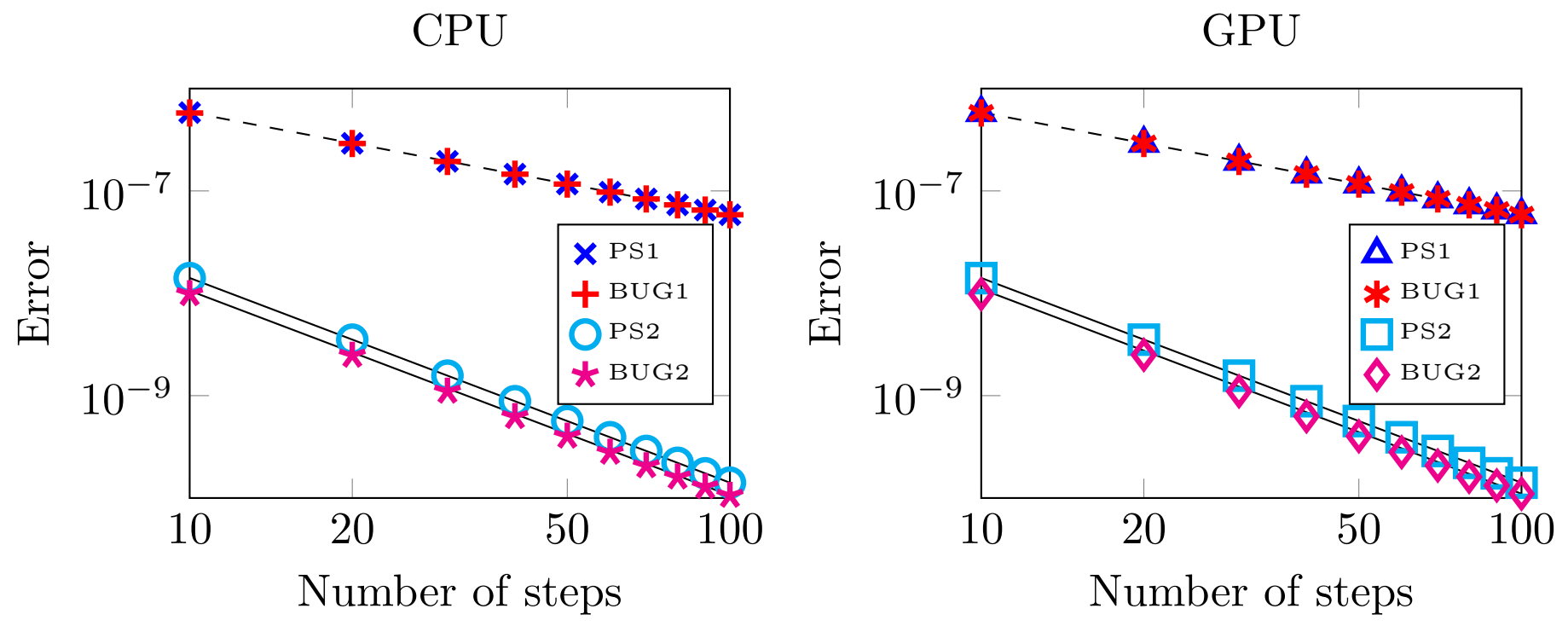
We start by assessing the order of convergence of the integrators to verify the
correctness of their implementation on both CPU and GPU architectures.
We consider four schemes: first- and second-order PS (labeled PS1 and PS2,
respectively), and first- and 
second-order BUG (labeled BUG1 and BUG2, respectively).
The tests are performed for the linear Landau damping problem~\eqref{eq:ll}
with $16^3\times 16^3$
degrees of freedom, up to the final time $T=0.1$.
Note that we consider a relatively small-sized problem in this experiment, as
our primary objective is to verify the correctness of the implementation.
Moreover, since the error is evaluated by constructing the
full distribution function, the number of degrees of freedom must be limited
(in subsequent experiments, we will consider substantially larger numbers of
degrees of freedom).
The number of time steps ranges from 10 to 100, while the rank is fixed to 10
for all simulations.
The error is measured using the relative infinity norm with respect to a reference
solution computed with PS2 using 3000 time steps. The results are reported in
Figure~\ref{fig:orders}.
As expected, all integrators exhibit the theoretically predicted orders of
convergence, confirming the correctness of their implementations.
Notice in addition that at first order, the PS and BUG schemes yield comparable errors.
A similar agreement is observed at second order, although BUG achieves slightly smaller errors.
Finally, the comparison between CPU and GPU implementations shows no appreciable
qualitative or quantitative differences.
In particular, the solutions obtained on the two architectures agree to nearly 
machine precision, with discrepancies appearing only around the thirteenth
digit.

\subsection{Comparison PS and BUG}
We now compare the proposed PS approach, considering adaptivity in terms of the
distribution function (\texttt{df} approach), with the BUG adaptive technique, as these are the two
most natural procedures to relate due to their similarities.
We consider both first-order and second-order integrators.
The tests are carried out for the two-stream instability problem~\eqref{eq:ts}
on an $80^3\times 80^3$ grid, using a time step size of $\tau=0.04$ and a final
time of $T=40$. For the PS integrators, the absolute tolerances for rank increase
and reduction are set to $\epsilon^\text{abs}_\text{incr}=10^{-4}$ and
$\epsilon^\text{abs}_\text{decr}=10^{-6}$, respectively. The relative tolerances
are set to $1/10$ of the absolute ones.
For the BUG adaptive methods, the needed tolerance is set to $10^{-4}$.
In all cases, the initial rank is set to 25, the minimum rank to 10, and the
maximum rank to 225.
Note that the maximum rank allowed for BUG is chosen to match that of PS.
This means that in practice for the former we have to deal with objects of column size
450 at first order and 675 at second order,
since the BUG methods require to build intermediate basis of size of $2r$ and $3r$,
respectively.
We also highlight that the same experiment with a grid size
of $90^3\times 90^3$ cannot be executed using the second-order BUG integrator
due to excessive memory requirements. Thus, the examples
presented in this section are essentially the largest we can consider on the available
hardware when employing BUG2. In contrast, the PS approach allows us to
perform considerably larger-scale experiments, as demonstrated in the following
sections.

\subsubsection{Wall-clock times for first-order integrators}\label{sec:fvsbug}
\begin{table}[htb]
  \centering
  \small
  \renewcommand{\arraystretch}{0.8}
  \begin{tabular}{l|cc||cc}
    & \multicolumn{2}{c}{\textit{BUG1}} & \multicolumn{2}{c}{\textit{PS1}} \\
    & GPU & CPU & GPU & CPU\\
    Initialization & $7.22\cdot 10^{-1}$ & $8.41\cdot 10^{-1}$ & $7.16\cdot 10^{-1}$ & $8.78\cdot 10^{-1}$ \\
    Time integration & $4.52\cdot 10^{2}$ & $8.55\cdot 10^{3}$ & $3.73\cdot 10^{2}$ & $6.86\cdot 10^{3}$ \\
    Rank choice & $1.17\cdot 10^{1}$ & $6.52\cdot 10^{1}$ & $7.77\cdot 10^{0}$ & $2.48\cdot 10^{1}$ \\\hhline{=|==||==}
    Total WC time & $4.64\cdot 10^2$ & $8.62\cdot 10^3$ & $3.81\cdot 10^2$ & $6.88\cdot 10^3$ \\
  \end{tabular}
  \caption{Summary of wall-clock times (in seconds) for the simulations using
  first-order adaptive BUG and PS \texttt{df}, see Section~\ref{sec:fvsbug}.}
  \label{tab:bdtimes_summ}
\end{table}
\begin{table}[htb]
  \centering
  \small
  \renewcommand{\arraystretch}{0.8}
  \begin{tabular}{l|c||l|c}
    \multicolumn{2}{c}{\textit{GPU}} & \multicolumn{2}{c}{\textit{CPU}} \\
                         & Time/step &                      & Time/step \\\hline
    V aug + QR           & $1.43\cdot 10^{-1}$ & X aug + QR           & $1.50\cdot 10^0$ \\
    $\quad \hookrightarrow$ \textit{QR part}           & $1.42\cdot 10^{-1}$ & $\quad \hookrightarrow$ \textit{QR part}           & $1.49\cdot 10^0$ \\
    X aug + QR           & $1.42\cdot 10^{-1}$ & V aug + QR           & $1.49\cdot 10^0$ \\
    $\quad \hookrightarrow$ \textit{QR part}           & $1.41\cdot 10^{-1}$ & $\quad \hookrightarrow$ \textit{QR part}           & $1.48\cdot 10^0$ \\
    L step               & $4.26\cdot 10^{-2}$ & L step               & $1.22\cdot 10^0$ \\
    K step               & $4.16\cdot 10^{-2}$ & K step               & $1.22\cdot 10^0$  \\
    D coeffs aug         & $3.05\cdot 10^{-2}$ & C coeffs aug         & $9.91\cdot 10^{-1}$ \\
    C coeffs aug         & $1.84\cdot 10^{-2}$ & D coeffs aug         & $9.90\cdot 10^{-1}$ \\
    S step               & $1.29\cdot 10^{-2}$ & D coeffs             & $4.91\cdot 10^{-1}$ \\
    SVD dec.             & $9.26\cdot 10^{-3}$ & C coeffs             & $4.91\cdot 10^{-1}$ \\
    D coeffs             & $8.47\cdot 10^{-3}$ & Hat matrices (for S) & $1.42\cdot 10^{-1}$ \\
    C coeffs             & $7.72\cdot 10^{-3}$ & SVD decomposition    & $4.37\cdot 10^{-2}$ \\
    Truncate and updates & $2.38\cdot 10^{-3}$ & Truncate and updates & $2.15\cdot 10^{-2}$ \\
    Hat matrices (for S) & $2.12\cdot 10^{-3}$ & S step               & $6.47\cdot 10^{-3}$ \\
    El.~field            & $1.09\cdot 10^{-3}$ & El.~field            & $5.17\cdot 10^{-3}$ \\
  \end{tabular}\\[2ex]
  \begin{tabular}{l|ccc||l|ccc}
    \multicolumn{4}{c}{\textit{GPU}} & \multicolumn{4}{c}{\textit{CPU}} \\
                  & Tot.~time &Calls& Time/call &                        & Tot.~time &Calls& Time/call \\\hline
    L step + QR   & $1.75\cdot 10^2$ &1205&    $1.45\cdot 10^{-1}$ &L step + QR   & $2.76\cdot 10^3$ & 1211 & $2.28\cdot 10^0$\\
    $\quad \hookrightarrow$ \textit{QR part}   & $1.11\cdot 10^2$ &1205&    $9.19\cdot 10^{-2}$ &$\quad \hookrightarrow$ \textit{QR part}   & $9.42\cdot 10^2$ & 1211 & $7.78\cdot 10^{-1}$\\
    K step + QR   & $1.73\cdot 10^2$ &1205&    $1.44\cdot 10^{-1}$ &K step + QR   & $2.75\cdot 10^3$ & 1211 & $2.27\cdot 10^0$\\
    $\quad \hookrightarrow$ \textit{QR part}   & $1.10\cdot 10^2$ &1205&    $9.16\cdot 10^{-2}$ &$\quad \hookrightarrow$ \textit{QR part}   & $9.43\cdot 10^2$ & 1211 & $7.78\cdot 10^{-1}$\\
    D coeff.      & $1.18\cdot 10^1$ &1205&    $9.79\cdot 10^{-3}$ &D coeff.      & $6.65\cdot 10^2$ & 1211 & $5.49\cdot 10^{-1}$\\
    C coeff.      & $1.08\cdot 10^1$ &1205&    $8.96\cdot 10^{-3}$ &C coeff.      & $6.64\cdot 10^2$ & 1211 & $5.49\cdot 10^{-1}$\\
    SVD dec.      & $6.18\cdot 10^0$ &1205&    $5.13\cdot 10^{-3}$ &SVD dec.      & $1.67\cdot 10^1$ & 1211 & $1.38\cdot 10^{-2}$\\
    Reject (inc.) & $1.53\cdot 10^0$ &205&     $7.48\cdot 10^{-3}$ &Reject (inc.) & $8.07\cdot 10^0$ & 211 & $3.83\cdot 10^{-2}$\\
    El.~field     & $1.51\cdot 10^0$ &1205&    $1.25\cdot 10^{-3}$ &El.~field     & $6.61\cdot 10^0$ & 1211 & $5.46\cdot 10^{-3}$\\
    S step        & $4.89\cdot 10^{-1}$ &1205& $4.06\cdot 10^{-4}$ &S step        & $4.48\cdot 10^0$ & 1211 & $3.70\cdot 10^{-3}$\\
    Accept (dec.) & $1.30\cdot 10^{-3}$ &5&    $2.60\cdot 10^{-4}$ &Accept (main.) & $9.03\cdot 10^{-4}$ & 995 & $9.08\cdot 10^{-7}$\\
    Accept (main.) & $5.28\cdot 10^{-4}$ &995&  $5.31\cdot 10^{-7}$ &Accept (dec.) & $1.25\cdot 10^{-4}$ & 5& $2.49\cdot 10^{-5}$\\
  \end{tabular}
  \caption{Detailed timings (in seconds) for the main parts of
  adaptive 
  first-order BUG (top)
  and PS \texttt{df} (bottom), see Section~\ref{sec:fvsbug}.}
  \label{tab:bdtimes_new}
\end{table}
We start by comparing the wall-clock times (in seconds) for the
first-order integrators PS1 and BUG1, executed on both CPU and GPU architectures.
A summary of the overall computational costs is provided in Table~\ref{tab:bdtimes_summ},
while Table~\ref{tab:bdtimes_new} reports more detailed results for the most
computationally relevant parts of the algorithms.

As highlighted in Table~\ref{tab:bdtimes_summ}, the initialization phase
accounts for a negligible fraction of the total simulation time.
This is expected, since this phase only includes setting up the initial condition,
allocating the required memory, initializing the BLAS and FFTW plans, and, when applicable,
transferring the initial condition from the CPU to the GPU. Notably, this transfer
represents the only communication between the CPU and GPU, as all subsequent
parts of the procedures are executed entirely on the selected hardware architecture.
The summary also stresses the computational advantage of GPUs over CPUs for
both PS and BUG in the time-integration and rank-selection stages.
As expected, these parts are slightly more computationally expensive for BUG,
which needs to operate on quantities of size $2r$. In terms of total wall-clock time,
PS is approximately 1.3 times faster than BUG on both CPU and GPU. Moreover,
the GPU achieves a speedup of approximately 18 times over the CPU for
both integrators, demonstrating the substantial benefits of GPU hardware.
Note also that all simulations exhibit comparable rank evolution, as illustrated
in Figure~\ref{fig:rank_f_vs_BUG} (left). During the linear phase, the ranks grow
at the same rate for all configurations. As the system enters the nonlinear
phase of the two-stream instability, the ranks show small variations in the
growth rate, while still being comparable. Clearly, some differences between
PS and BUG are expected, since, despite their similarities, the two methods do not
employ exactly the same adaptivity strategy. Small discrepancies between CPU and
GPU simulations are also expected, especially for PS, due to the inherent randomness introduced when
augmenting the approximation space with random columns
(then orthonormalized as part of the rank-increase procedure).

Focusing now on Table~\ref{tab:bdtimes_new}, we observe that 
for the BUG integrator (top part of the table), the most computationally expensive components are the QR factorizations in the augmentation steps for $X$ and $V$.
Together, they account for approximately two-thirds of the overall run time on the GPU. This indicates that the implementation operates relatively close to what we can reasonably expect from an optimal implementation. The next most expensive components are the K and L steps, followed by the updates of the augmented C and D coefficients.
For PS (bottom part of the table), we report both the total wall-clock times
and the average times per function call. In this case, the normalization is
performed with respect to the total number of calls, which depends on
the outcome of the adaptive time-stepping procedure. Specifically,
for the GPU simulation, we have 205 rejected steps,
995 accepted steps that preserve the rank, and 5 accepted steps that reduce the
rank, yielding a total of 1205 actual numerical steps.
For the CPU simulation, the corresponding numbers are 211 rejected steps, 995 accepted steps that
preserve the rank, and 5 accepted steps that reduce the rank,
resulting in a total of 1211 actual numerical steps.
As expected, the most computationally expensive stages in PS are the L and
K steps. The QR factorizations account for approximately 60\% of the overall run time on the GPU. The K and L steps are followed by the updates of
the C and D coefficients. Rejected time steps
are more computationally expensive than accepted ones, whereas the S step
contributes negligibly to the overall computational cost. We also observe that
the costs of the K and L steps are comparable for BUG and PS. This is consistent with the similar rank evolution observed for the two integrators
(both steps operate on quantities with $r$ columns).

Finally, the speedups achieved by the GPU relative to the CPU
across the different parts of the procedure are as expected. Therefore, all subsequent simulations will
be performed exclusively on the GPU.

\subsubsection{Second-order integrators comparison}\label{sec:fvsbug2}
\begin{figure}[htb]
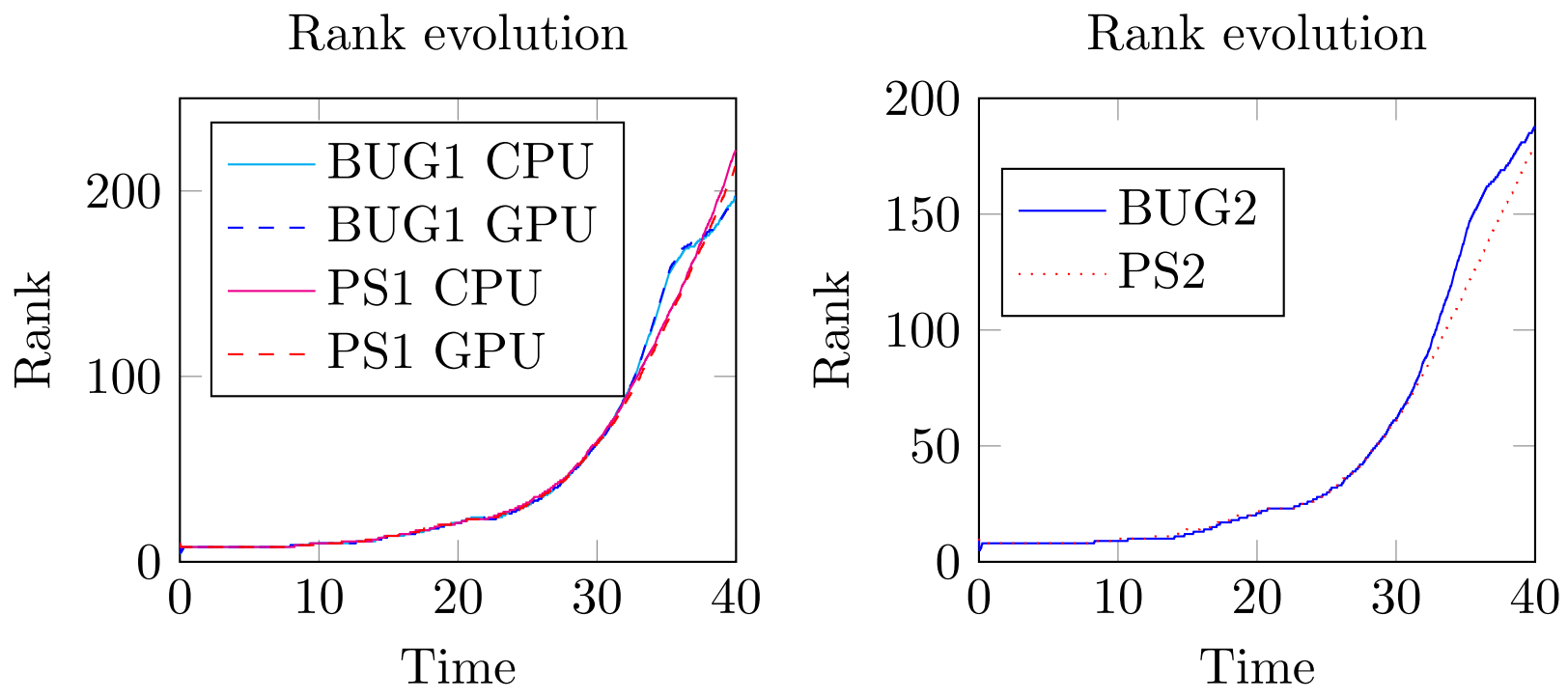

  \centering
  \input{plot_ranks_cpugpu.tex}
  \input{plot_rank_f_vs_bug.tex}
  \caption{Rank evolution for adaptive BUG and PS \texttt{df}.
  Left: first order on CPU and GPU (see Section~\ref{sec:fvsbug}).
  Right: second order on GPU (see Section~\ref{sec:fvsbug2}).}
  \label{fig:rank_f_vs_BUG}
\end{figure}

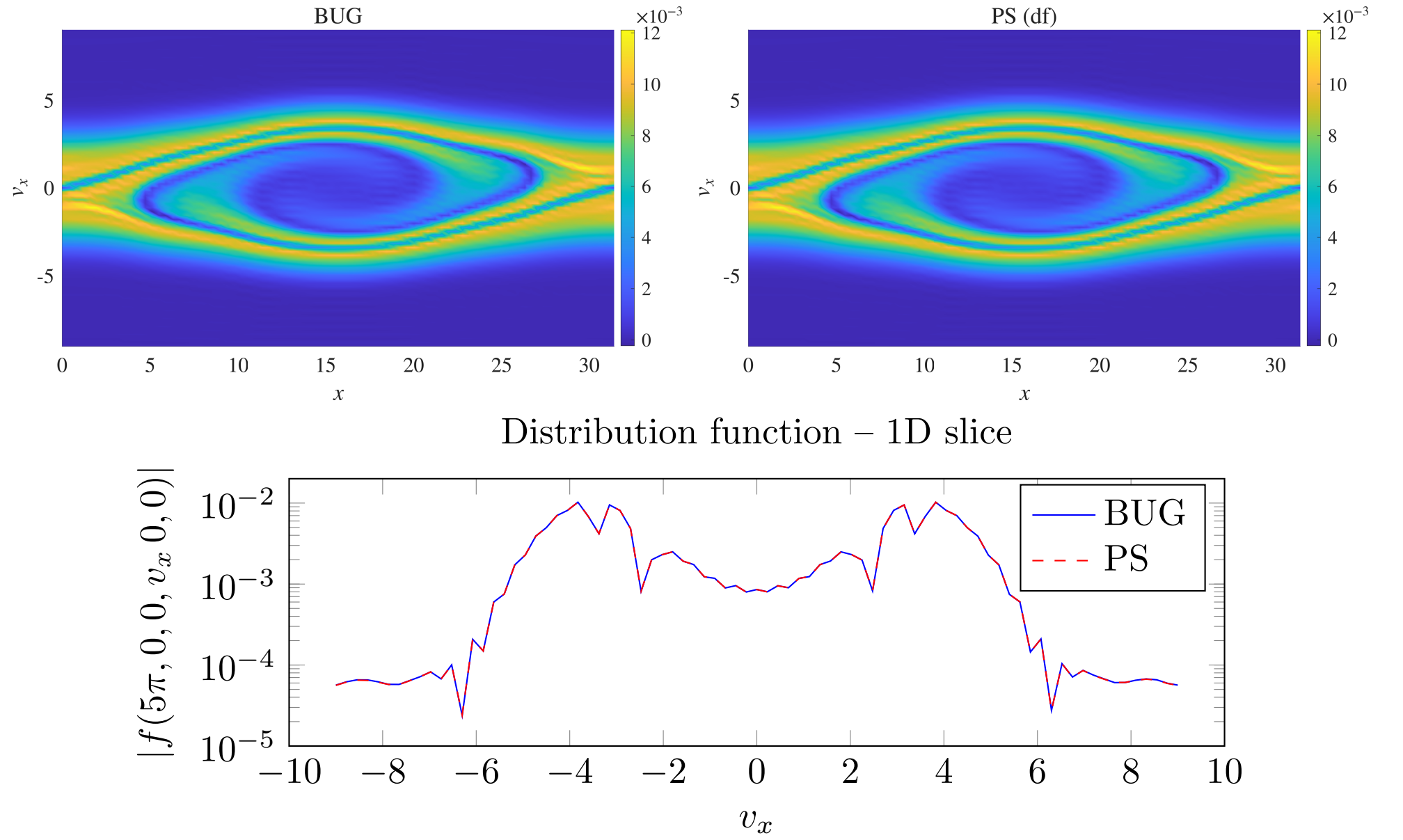
\begin{figure}[htb]
  \centering
  \includegraphics[scale=0.225]{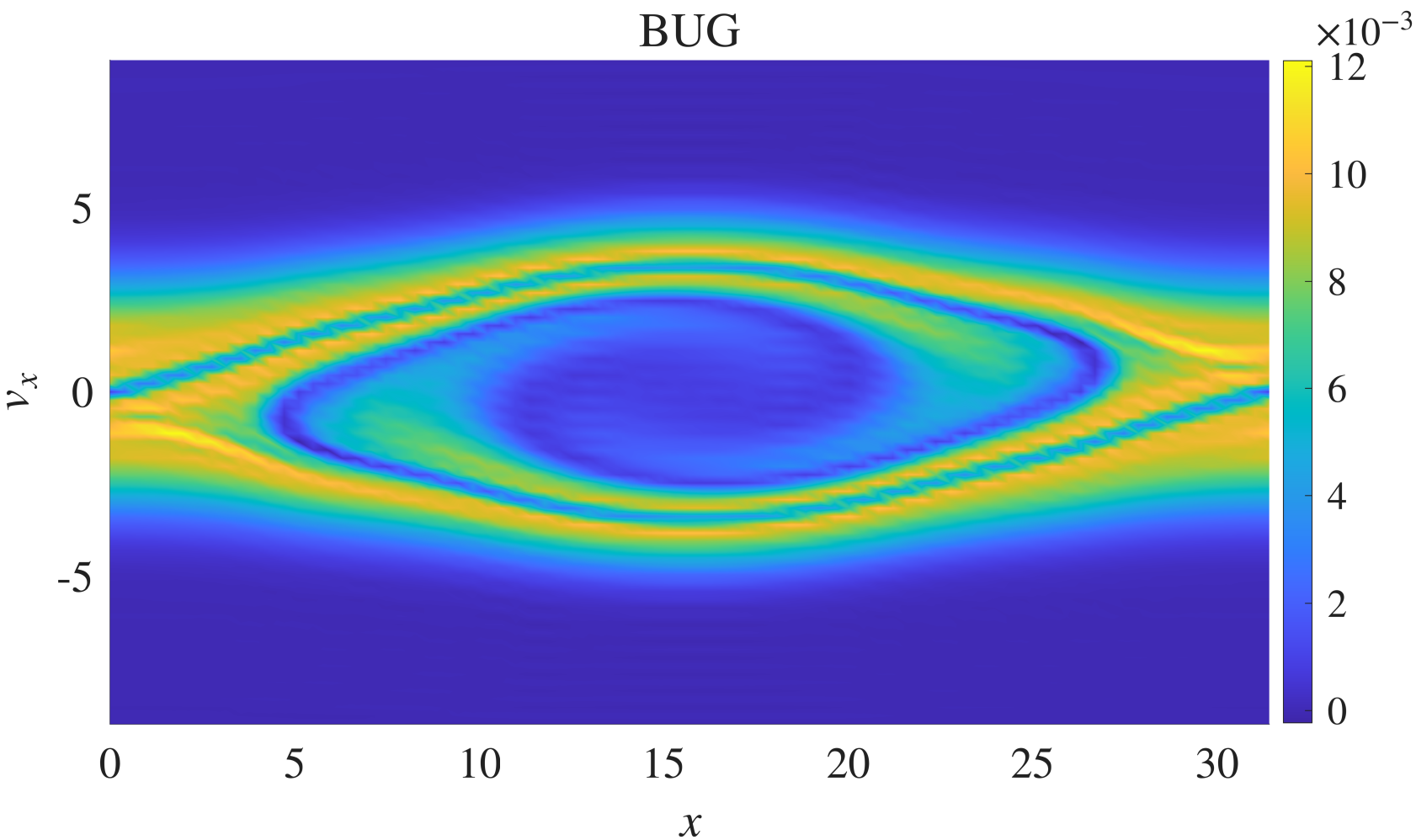}
  \includegraphics[scale=0.225]{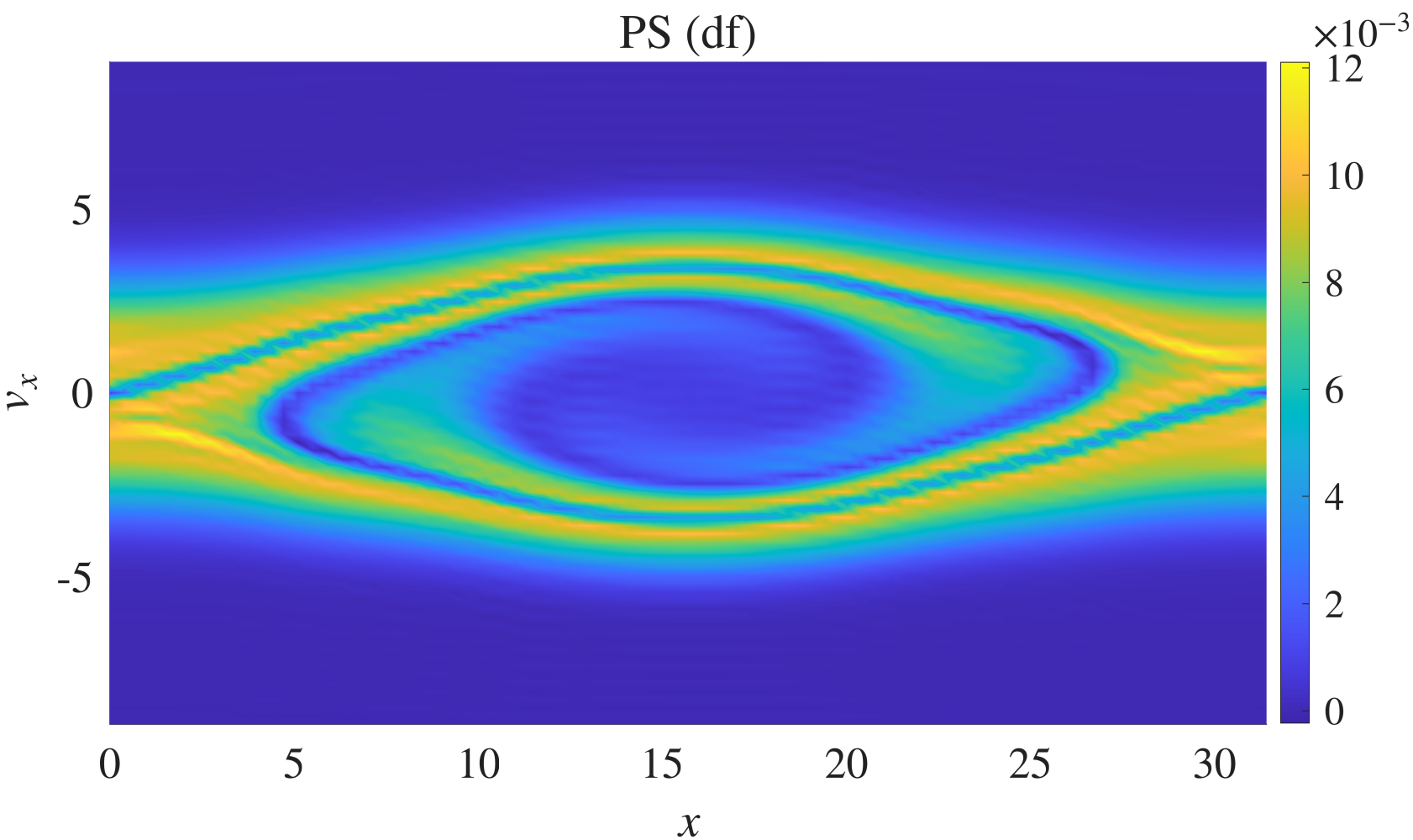}
  \\
%
%
%
\begin{tikzpicture}

\begin{axis}[%
width=3.5in,
height=1in,
at={(1.619in,0.822in)},
scale only axis,
xmin=-10,
xmax=10,
xlabel={$v_x$},
ymode=log,
ymin=1e-05,
ymax=2e-02,
yminorticks=true,
ylabel={$|f(5\pi,0,0,v_x\,0,0)|$},
axis background/.style={fill=white},
title={Distribution function -- 1D slice},
legend style={legend cell align=left, align=left}
]
\addplot [color=blue]
  table[row sep=crcr]{%
-9	5.647021115161e-05\\
-8.775	6.2030743934783e-05\\
-8.55	6.56128377643333e-05\\
-8.325	6.55152453056186e-05\\
-8.1	6.20609229106507e-05\\
-7.875	5.77152466832158e-05\\
-7.65	5.7610857810637e-05\\
-7.425	6.41913485369724e-05\\
-7.2	7.16891618226762e-05\\
-6.975	8.23779224954836e-05\\
-6.75	6.73766070845843e-05\\
-6.525	0.000100439353236463\\
-6.3	2.38161075067478e-05\\
-6.075	0.000207461622511615\\
-5.85	0.000149314275652508\\
-5.625	0.000600045109557177\\
-5.4	0.000753259388236967\\
-5.175	0.00172741892845356\\
-4.95	0.00228665488356493\\
-4.725	0.00391575414182707\\
-4.5	0.00497561610217612\\
-4.275	0.00701263086833471\\
-4.05	0.00810959553096537\\
-3.825	0.0102507650358637\\
-3.6	0.00679245007986868\\
-3.375	0.00418217345186032\\
-3.15	0.00949495485380588\\
-2.925	0.00810551392731971\\
-2.7	0.00484646564104407\\
-2.475	0.000814294678016049\\
-2.25	0.00199336783490067\\
-2.025	0.00230479882799127\\
-1.8	0.00250159697165177\\
-1.575	0.00191448294515547\\
-1.35	0.00174076820963383\\
-1.125	0.00122913870933006\\
-0.9	0.00117736566640491\\
-0.675000000000001	0.000895088979028376\\
-0.449999999999999	0.000956324194665859\\
-0.225	0.00079787518249375\\
0	0.000855997651554282\\
0.225	0.000802577454910952\\
0.449999999999999	0.000953236354035377\\
0.675000000000001	0.000900367543003046\\
0.9	0.00117369388259317\\
1.125	0.00123660988773241\\
1.35	0.00173827119257279\\
1.575	0.00192491218897943\\
1.8	0.00249881756980049\\
2.025	0.00231236156927239\\
2.25	0.0019789024353817\\
2.475	0.000823803502846939\\
2.7	0.0048593971399456\\
2.925	0.00812137913580025\\
3.15	0.00948048963969289\\
3.375	0.00417338352853473\\
3.6	0.0068265373449466\\
3.825	0.0102473916480988\\
4.05	0.00810044974098493\\
4.275	0.00700529249641728\\
4.5	0.00496791809692414\\
4.725	0.00391014869232457\\
4.95	0.0022787424274576\\
5.175	0.00172490592227439\\
5.4	0.000747042084475469\\
5.625	0.000601530948487168\\
5.85	0.000144965602316489\\
6.075	0.000210015779057754\\
6.3	2.76067295755571e-05\\
6.525	0.000103969764815339\\
6.75	7.10934054390594e-05\\
6.975	8.58300299332755e-05\\
7.2	7.49707180940209e-05\\
7.425	6.74209094799191e-05\\
7.65	6.06731435069186e-05\\
7.875	6.07769658536421e-05\\
8.1	6.49004008281264e-05\\
8.325	6.74042272330921e-05\\
8.55	6.56879140400517e-05\\
8.775	5.96443461186994e-05\\
9	5.647021115161e-05\\
};
\addlegendentry{BUG}

\addplot [color=red, dashed]
  table[row sep=crcr]{%
-9	5.61148069098316e-05\\
-8.775	6.15812946058925e-05\\
-8.55	6.5096497674415e-05\\
-8.325	6.46328929732109e-05\\
-8.1	6.1426913319696e-05\\
-7.875	5.75552658837017e-05\\
-7.65	5.73990305882485e-05\\
-7.425	6.39380154793424e-05\\
-7.2	7.16895332492627e-05\\
-6.975	8.22357218776496e-05\\
-6.75	6.7328462461958e-05\\
-6.525	9.99560394230696e-05\\
-6.3	2.41774793989171e-05\\
-6.075	0.000208603757899584\\
-5.85	0.00014917831892808\\
-5.625	0.000600488691492267\\
-5.4	0.000754834029822667\\
-5.175	0.00172793349682322\\
-4.95	0.00228902871042452\\
-4.725	0.00391663775642054\\
-4.5	0.00497758165420609\\
-4.275	0.00701283639048726\\
-4.05	0.00811056887805989\\
-3.825	0.0102511980970133\\
-3.6	0.00678944186361479\\
-3.375	0.00418514738502707\\
-3.15	0.00949637206199892\\
-2.925	0.00809964622717593\\
-2.7	0.00483846254620651\\
-2.475	0.000816447267255403\\
-2.25	0.00199647436062302\\
-2.025	0.0023056518923807\\
-1.8	0.00250243548523002\\
-1.575	0.00191607388724761\\
-1.35	0.0017424170853223\\
-1.125	0.00123324607093164\\
-0.9	0.00117730158614233\\
-0.675000000000001	0.000896038463060782\\
-0.449999999999999	0.000956932413821512\\
-0.225	0.000795473863510359\\
0	0.000855846827962251\\
0.225	0.000800207687225346\\
0.449999999999999	0.000951345289003148\\
0.675000000000001	0.000900716865400838\\
0.9	0.00117498699783873\\
1.125	0.0012394846814884\\
1.35	0.00173546067510414\\
1.575	0.0019250466709367\\
1.8	0.00249654035354392\\
2.025	0.00231241988174308\\
2.25	0.00197592670803412\\
2.475	0.000831631079152927\\
2.7	0.00486757419555052\\
2.925	0.00812532290751599\\
3.15	0.00947841098459777\\
3.375	0.00416068697220331\\
3.6	0.00683003358895266\\
3.825	0.0102470804892503\\
4.05	0.00810224505070391\\
4.275	0.00700456032957148\\
4.5	0.00496861470138836\\
4.725	0.00390892360656046\\
4.95	0.00227831286536867\\
5.175	0.00172484056810512\\
5.4	0.000748084859953072\\
5.625	0.000602637516174634\\
5.85	0.00014513055245442\\
6.075	0.0002097971045855\\
6.3	2.7347554625737e-05\\
6.525	0.000103511177924949\\
6.75	7.10383069885118e-05\\
6.975	8.54998790286882e-05\\
7.2	7.43404688304559e-05\\
7.425	6.68842181667005e-05\\
7.65	6.09124620594396e-05\\
7.875	6.09579037975008e-05\\
8.1	6.4548520915659e-05\\
8.325	6.69241030845173e-05\\
8.55	6.51513514392037e-05\\
8.775	5.93230554976391e-05\\
9	5.61148069098316e-05\\
};
\addlegendentry{PS}

\end{axis}

\end{tikzpicture}%
  \caption{Plot of 2D slices $|f(x,0,0,v_x,0,0)|$ (top) and 1D slice
  $|f(5\pi,0,0,v_x,0,0)|$ (bottom) at the final time $T=40$ obtained with adaptive
  second-order BUG and PS,
  see the experiment in Section~\ref{sec:fvsbug2}.}
  \label{fig:distrib_fBUG}
\end{figure}

We now repeat the previous experiment using second-order integrators. For brevity
of exposition, we report only the total wall-clock times (in seconds), which 
for BUG midpoint is $8.08\cdot 10^{2}$ while for PS $6.85\cdot 10^{2}$
(yielding once again a speedup of approximately 1.2).
As in the previous experiment, the two methods exhibit comparable rank
evolution, as shown in Figure~\ref{fig:rank_f_vs_BUG} (right). To further assess the
agreement between the two simulations, we compare the corresponding
distribution functions at the final time $T=40$ in terms of both a two-dimensional phase space slice
(namely, $|f(x,0,0,v_x,0,0)|$)
and a one-dimensional slice (namely, $|f(5\pi,0,0,v_x,0,0)|$),
as shown in Figure~\ref{fig:distrib_fBUG}.
The two methods yield similar results for all the slices,
confirming that they provide comparable approximations of the distribution
function. This agreement is particularly relevant because the rank-adaptivity
strategy is driven by information extracted from the distribution function
itself (see the following section for comparisons using adaptivity based on different physical quantities).

To summarize, BUG provides a viable alternative to PS, showing comparable
simulation results in terms of computational times, rank evolution, and
distribution function. Nevertheless, PS offers a clear advantage within our
large-scale GPU framework, primarily because of its lower memory requirements.
In fact, by further reducing them, PS enables simulations at higher ranks and
with more degrees of freedom, as we will demonstrate in the following sections.

\subsection{Adaptivity based on physical quantities}\label{sec:diffadapt}
\begin{table}[htb]
  \centering
  \begin{tabular}{l|c|c|c}
                & $\epsilon_{\text{incr}}=10^{-3}$ & $\epsilon_{\text{incr}}=10^{-4}$ & $\epsilon_{\text{incr}}=10^{-5}$ \\\hline
    Adaptivity \texttt{df} &  $1.91\cdot 10^{3}$  & ---  & ---  \\
    Adaptivity \texttt{ef}  & $1.31\cdot 10^{3}$  & $2.06\cdot 10^{3}$  & $2.90\cdot 10^{3}$  \\
    Adaptivity \texttt{ee}  & $7.12\cdot 10^{2}$  & $8.86\cdot 10^{2}$  & $1.17\cdot 10^{3}$  \\
  \end{tabular}
  \caption{Wall-clock times (in seconds) for the numerical experiments 
  using second-order PS with different
  tolerances and adaptivities, see Section~\ref{sec:diffadapt}.}
  \label{tab:wctimes_tol_new}
\end{table}
In this section, we investigate the effect of different adaptivity criteria and
tolerance levels on the accuracy and computational cost of the simulations.
In particular, as discussed previously in the paper, we consider adaptivity
driven by the electric field (\texttt{ef} approach) and the
electric energy (\texttt{ee} approach).
The test problem is the two-stream instability~\eqref{eq:ts}, with a
phase space discretization of $128^3\times 128^3$. All simulations are
performed using the second-order adaptive PS integrator on the GPU, with
an initial rank of 25, a minimum rank of 10, and a maximum rank of 100.
We consider $\tau=0.032$ and the final time $T=40$.
The absolute tolerance for rank decrease is fixed
at $\epsilon^\text{abs}_{\text{decr}}=10^{-8}$,
while three
different absolute tolerances $\epsilon^\text{abs}_{\text{incr}}$ are considered for
rank increase, namely $10^{-3}$, $10^{-4}$, and $10^{-5}$.
The relative tolerances are set to one tenth of the corresponding absolute ones.
For notational purposes, from now on we will simply denote $\epsilon^\text{abs}_{\text{incr}}$
as $\epsilon_{\text{incr}}$.
Table~\ref{tab:wctimes_tol_new} summarizes the total wall-clock times for all simulations.

\begin{figure}[htb]
  \centering
  \input{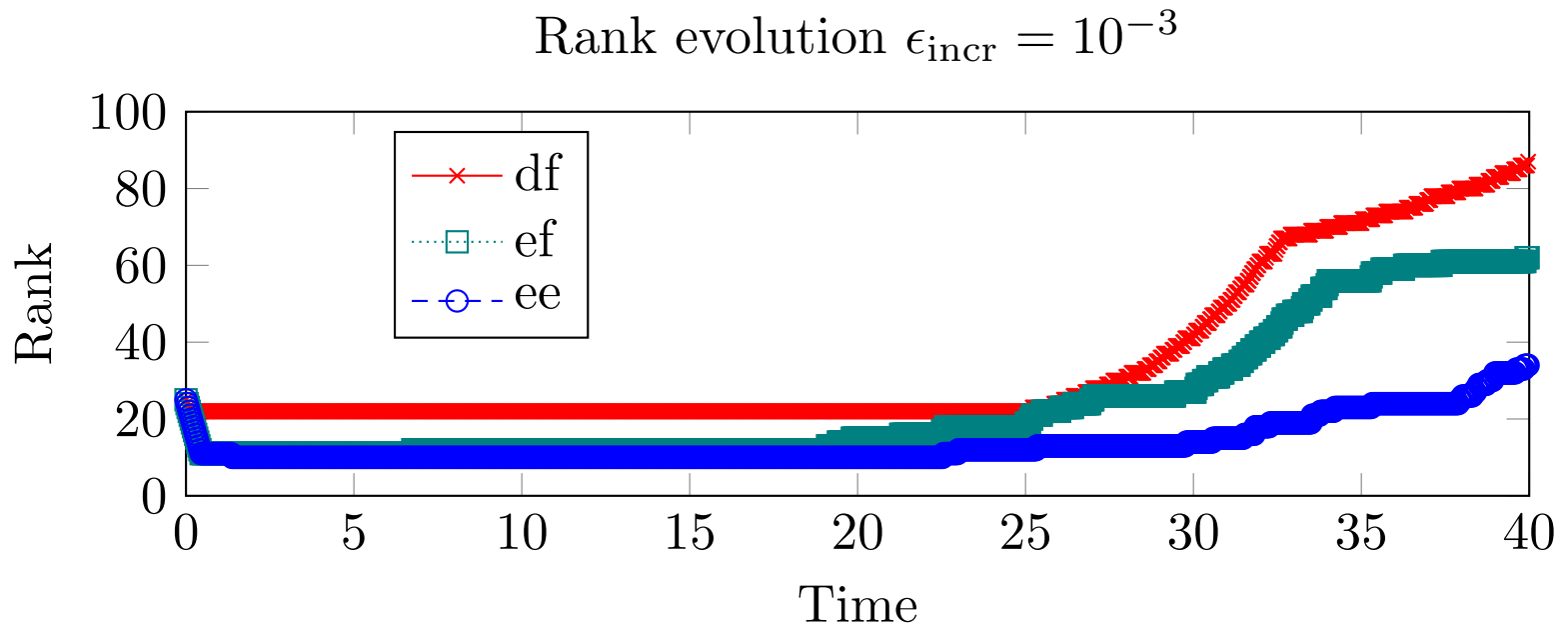}
  \caption{Rank evolution for second-order PS
  with $\epsilon_{\text{incr}}=10^{-3}$ and different adaptivity strategies,
  see Section~\ref{sec:diffadapt}.}
  \label{fig:ranks_1em3}
\end{figure}
\begin{figure}[htb]
  \centering
  \includegraphics[scale=0.225]{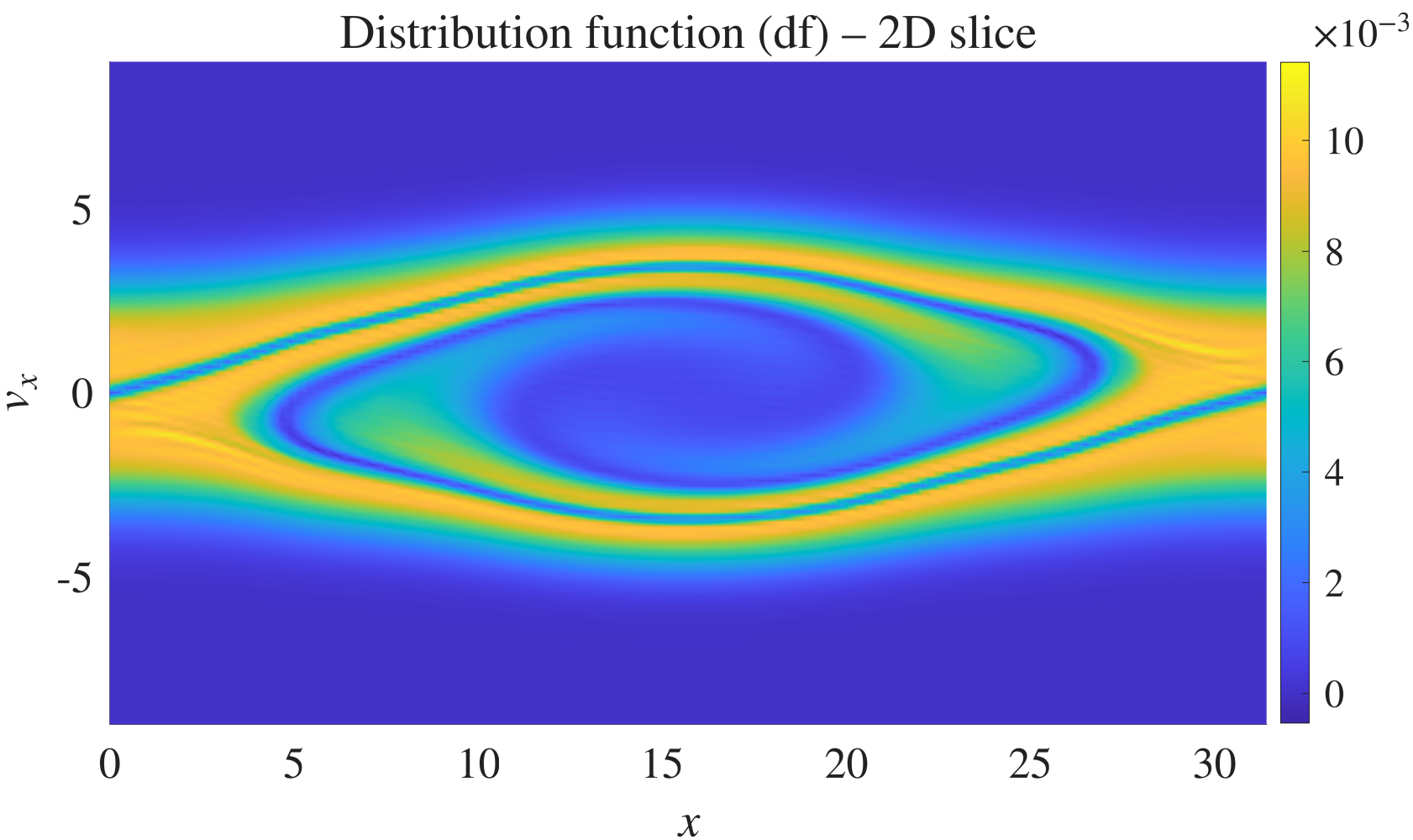}
  \includegraphics[scale=0.225]{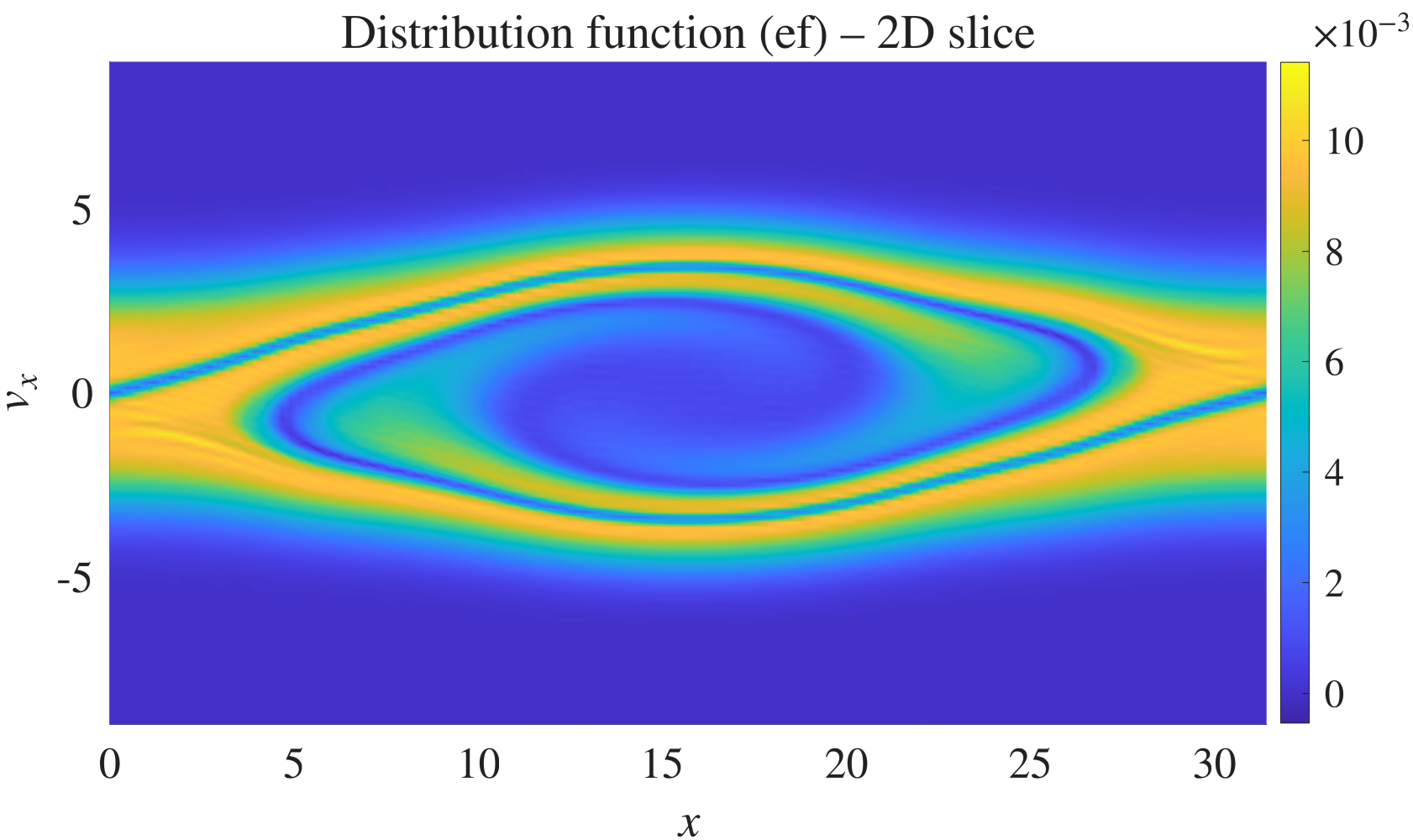}
  \\
  \includegraphics[scale=0.225]{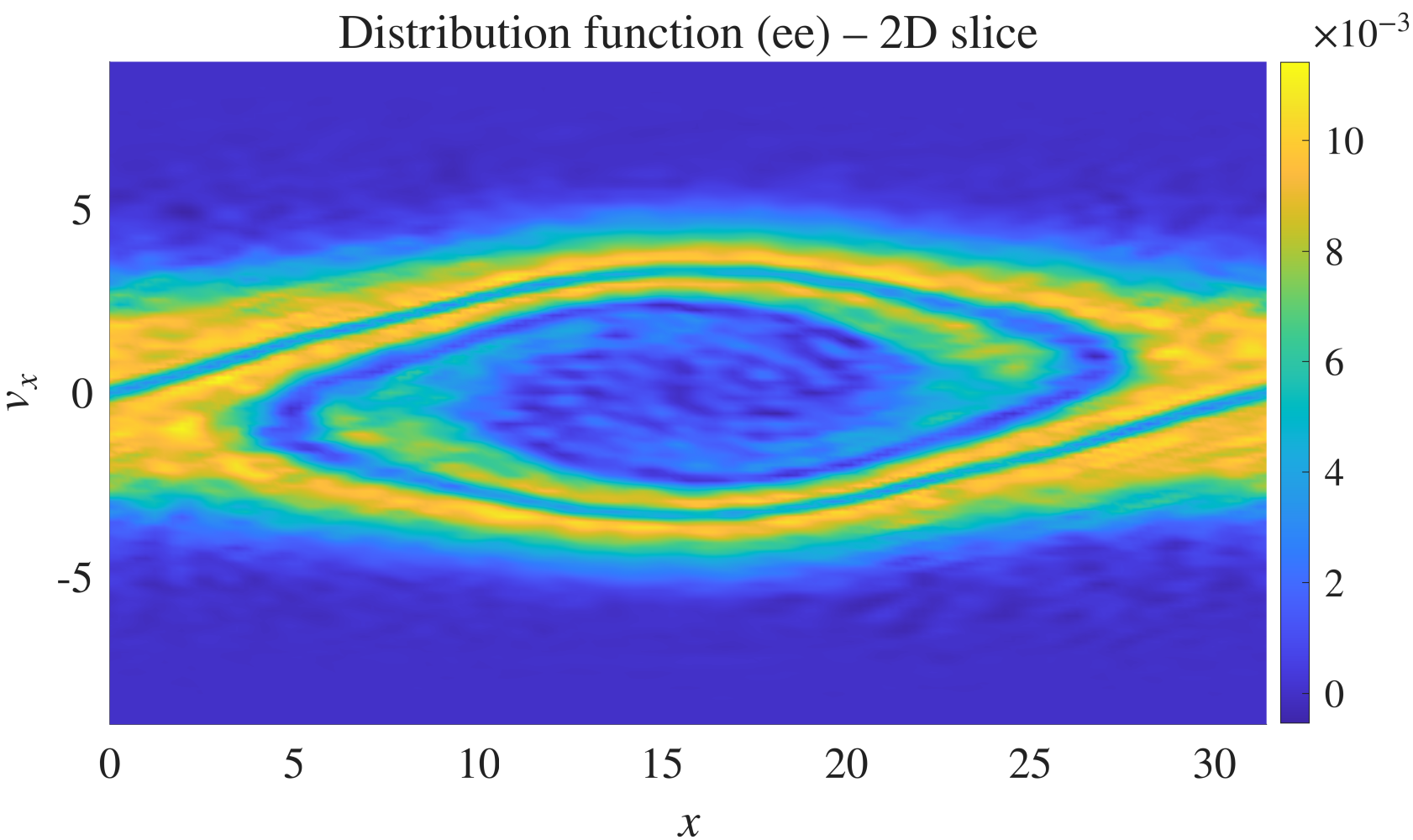}
  \input{plot_1dslice_f_ef_ee_1em3.tex}
  \caption{Plot of 2D slices $|f(x,0,0,v_x,0,0)|$ and 1D slice
  $|f(5\pi,0,0,v_x,0,0)|$ at the final time $T=40$ obtained with second-order PS and different
  adaptivity strategies ($\epsilon_{\text{incr}}=10^{-3}$),
  see Section~\ref{sec:diffadapt}.}
  \label{fig:distrib_1em3}
\end{figure}
\begin{figure}[htb]
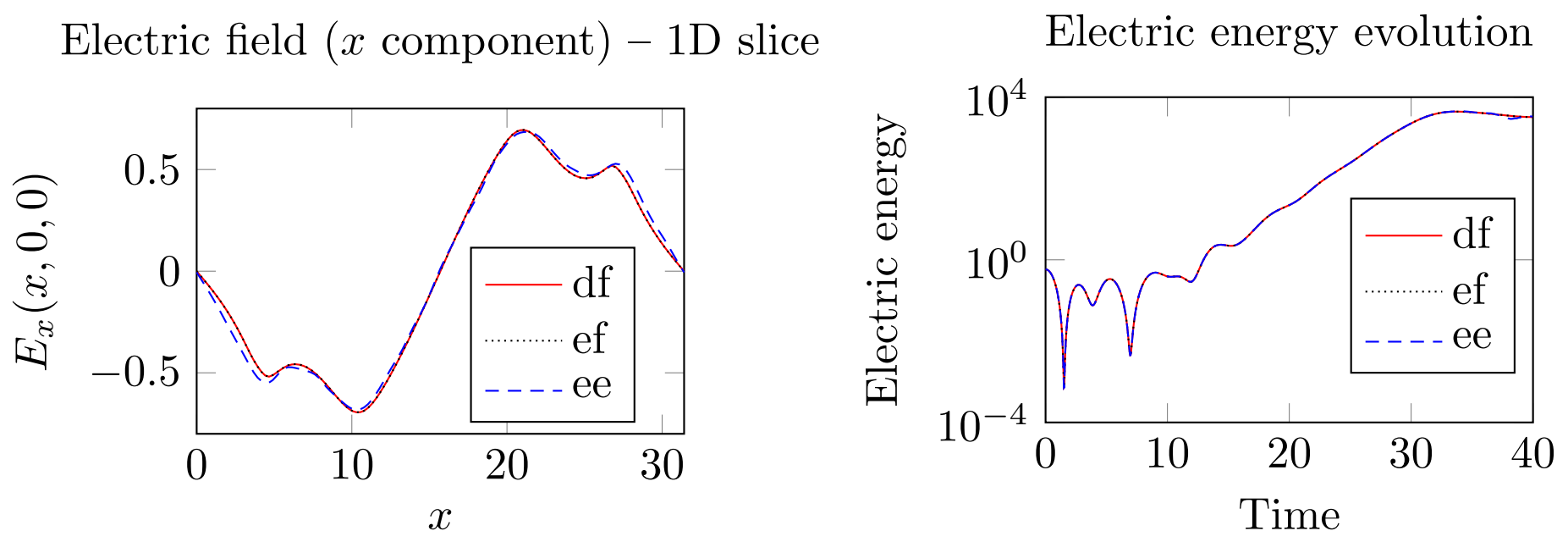

  \centering
  \input{plot_elfields_f_ef_ee_1em3.tex}
  \input{plot_elenergy_f_ef_ee_1em3.tex}
  \caption{Plot of $x$-component of the electric field at the final time $T=40$
  (left)
  and evolution of the electric energy (right)
  obtained with second-order PS, different
  adaptivity strategies, and $\epsilon_{\text{incr}}=10^{-3}$,
  see Section~\ref{sec:diffadapt}.}
  \label{fig:elfieldsenergies_1em3}
\end{figure}

We first discuss the results obtained with $\epsilon_{\text{incr}}=10^{-3}$.
For comparison, we also include the results obtained using the distribution function as the adaptivity criterion
(approach \texttt{df}).
The corresponding rank evolutions are shown in Figure~\ref{fig:ranks_1em3}.
As expected, driving adaptivity through physical quantities, namely the
electric field or the electric energy, leads to lower ranks than using
the distribution function. This reduction in rank directly translates into lower
computational and memory costs.
To assess the accuracy of the different adaptivity strategies, in 
Figure~\ref{fig:distrib_1em3} we compare two-dimensional slices of
the distribution function $|f(x,0,0,v_x,0,0)|$, together with the
corresponding one-dimensional slices $|f(5\pi,0,0,v_x,0,0)|$.
In addition, Figure~\ref{fig:elfieldsenergies_1em3} shows the $x$-component of
the electric field, evaluated at $y=z=0$ and plotted as a function of $x$,
together with the evolution of the electric energy.
Notice that we only plot the $x$-component of the electric field since it
is the dominant contribution (in fact, the other two components are
approximately two orders of magnitude smaller).
The results highlight clear differences between the adaptivity strategies.
When adaptivity is controlled by the electric energy, the overall evolution of
this quantity is adequately captured in both the linear and nonlinear phases of
the instability. However, this strategy leads to noticeable discrepancies in
the distribution function and some visible deviations in the electric field.
In contrast, controlling adaptivity through the electric field yields a more
accurate representation of the field itself and also captures the evolution of
the electric energy. Although some small differences remain in the
distribution function, they are less pronounced than those observed when using
the electric energy as the adaptivity criterion (see particularly the 1D slice).

\begin{figure}[htb]
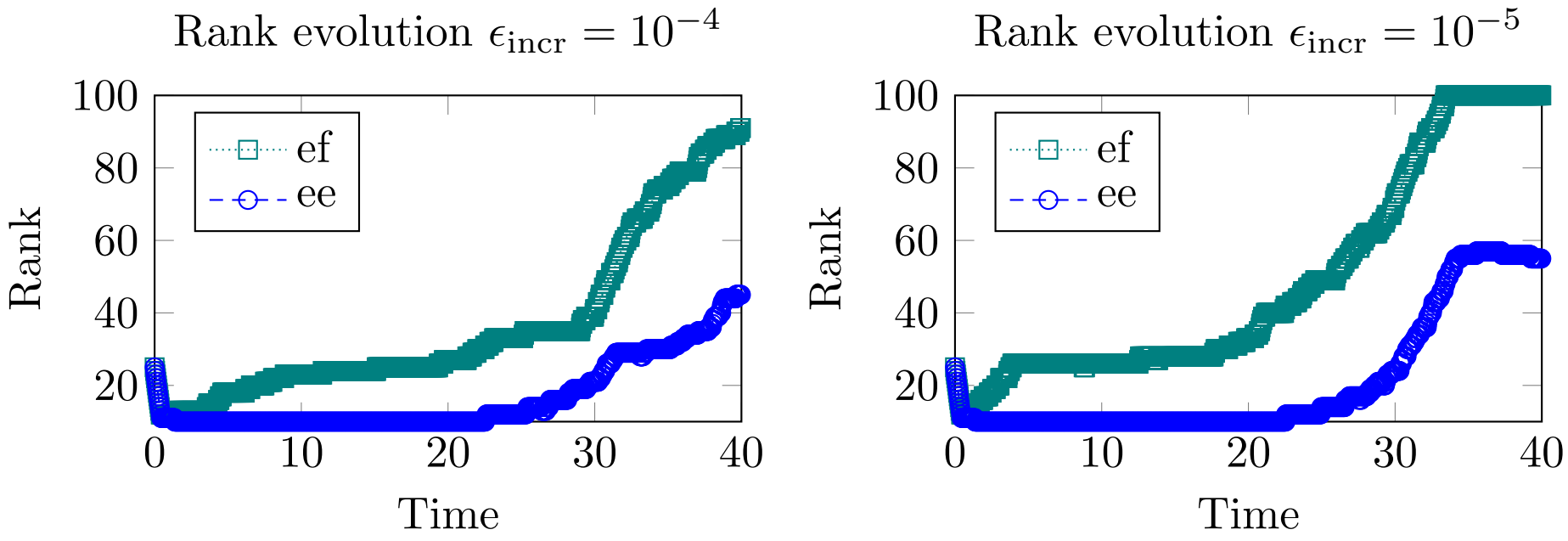

  \centering
  \input{plot_rank_1em4.tex}
  \input{plot_rank_1em5.tex}
  \caption{Rank evolution for second-order PS
  with $\epsilon_{\text{incr}}=10^{-4}$ (left) and
  $\epsilon_{\text{incr}}=10^{-5}$ (right) using different adaptivity strategies,
  see Section~\ref{sec:diffadapt}.}
  \label{fig:ranks_1em41em5}
\end{figure}
Finally, we repeat the experiment using increased tolerances of $10^{-4}$ and
$10^{-5}$. For brevity, we report only the corresponding rank evolutions in
Figure~\ref{fig:ranks_1em41em5}. As expected, tighter tolerances lead to
higher ranks, reflecting the increased accuracy requirements. In all instances,
the ranks obtained when controlling adaptivity through the electric energy
remain lower than those obtained with electric field-based adaptivity,
which is consistent with the physical behavior and the results observed for the rank-increase tolerance of $10^{-3}$.

\subsection{Benchmark examples}
We conclude the numerical examples by presenting the results of the proposed
rank-adaptive approach on a set of classical 6D benchmark experiments.
For each simulation, we report the electric-energy dynamics and
the rank evolution over time.
In all cases, we employ the presented second-order PS
integrator on the GPU, with adaptivity driven by the electric energy (approach \texttt{ee}),
which is the least computationally expensive. The tolerances reported here are absolute,
and the corresponding relative tolerances are always set to $1/10$ of those values.

\subsubsection{Linear Landau problem}\label{sec:ll}
\begin{figure}[htb]
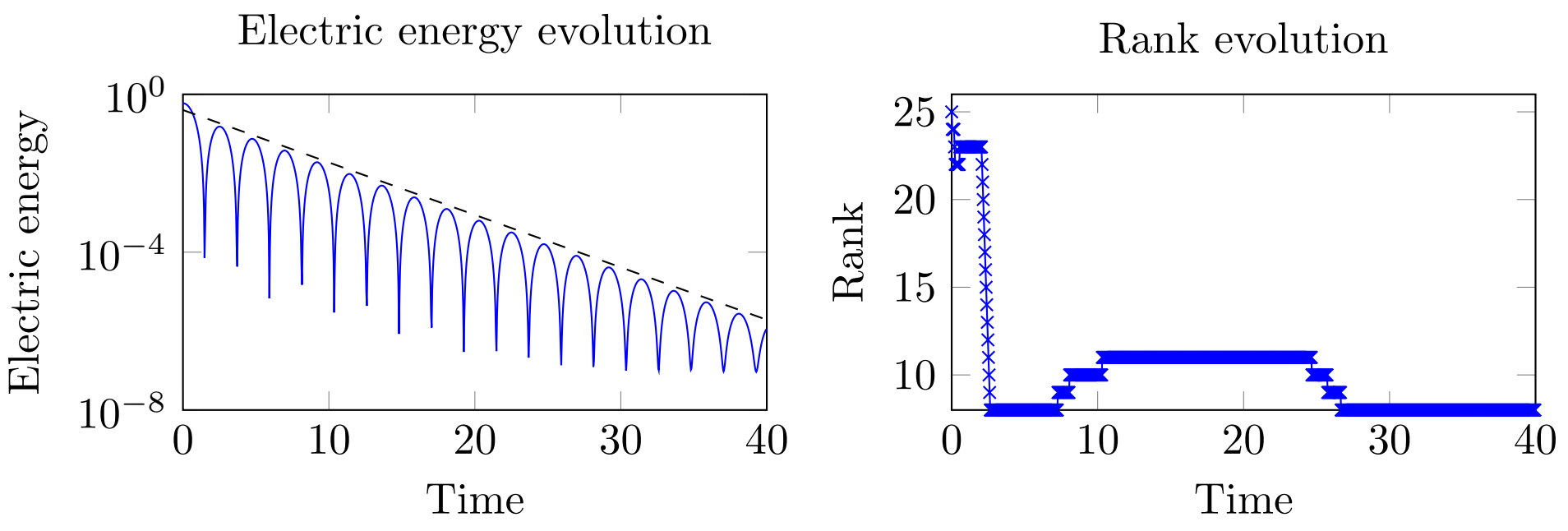

  \centering
  \input{plot_ee_ll.tex}
  \input{plot_rank_ll.tex}
  \caption{Results for the $64^3\times 256^3$ linear Landau benchmark example
  (see Section~\ref{sec:ll} for more details).
  Left: electric energy evolution (with reference line of slope $-0.306$).
  Right: rank evolution.}
  \label{fig:flagship_ll}
\end{figure}

The first simulation that we present is a classical linear Landau damping problem,
whose initial condition is detailed in equation~\eqref{eq:ll}. We set the
space and velocity degrees of freedom to $64^3\times 256^3$,
while the relevant tolerances for our algorithm are set to
$\epsilon_{\text{incr}}=10^{-6}$ and $\epsilon_{\text{decr}}=10^{-10}$.
Moreover, the initial, minimum, and maximum ranks are set to
the values 25, 8, and 25, respectively. The final time is $T=40$ and the time
step size is chosen as $\tau=0.04$.
The results are depicted in Figure~\ref{fig:flagship_ll}. As expected, the
electric energy decreases over time with a damping rate that matches the
theoretical estimate. In terms of rank adaptivity, we observe that in the
initial phase the rank decreases, reaching the minimum rank imposed and then
remaining between 8 and 11. This is expected, since in this
setting we are in a linear regime, hence high ranks are unlikely to occur.
The wall-clock time of the overall simulation is approximately 20 minutes. Note that
for the selected number of degrees of freedom, a full-rank simulation would
require approximately 35 TB of memory just to store the initial condition, making the
problem infeasible without a large supercomputing cluster.

\subsubsection{Nonlinear Landau problem}\label{sec:nll}
\begin{figure}[htb]
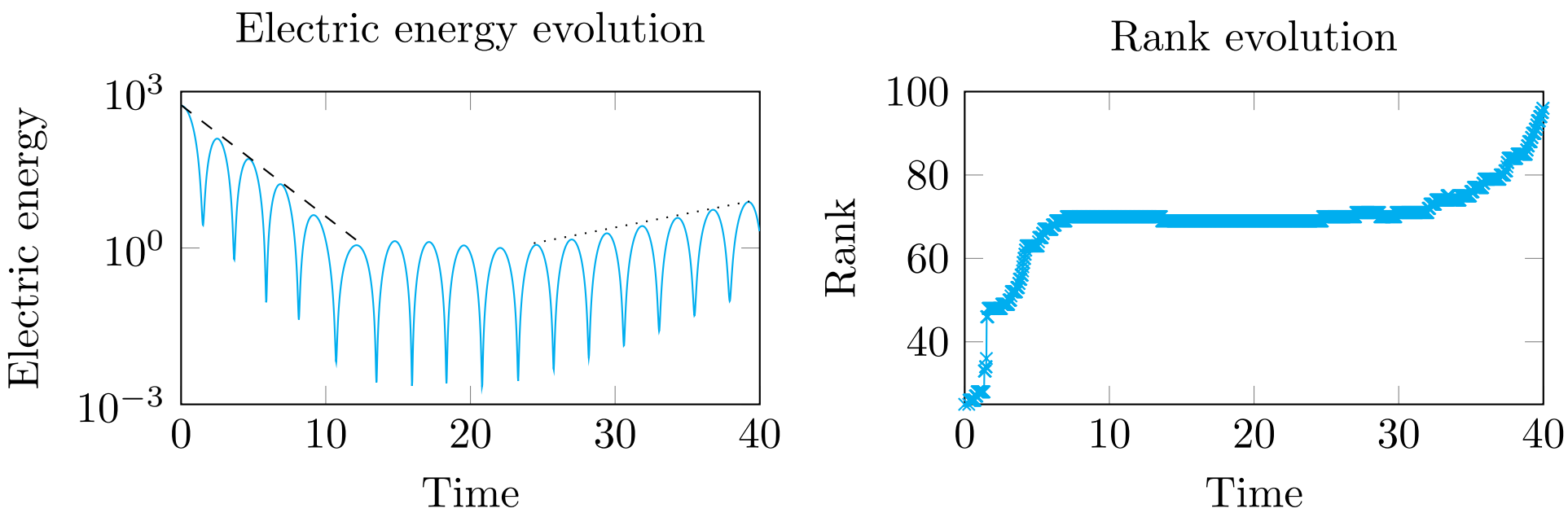

  \centering
  \input{plot_ee_nll.tex}
  \input{plot_rank_nll.tex}
  \caption{Results for the $64^3\times 128^3$ nonlinear Landau benchmark example
  (see Section~\ref{sec:nll} for more details).
  Left: electric energy evolution (with reference lines of slope $-0.492$ and $0.124$).
  Right: rank evolution.}
  \label{fig:flagship_nll}
\end{figure}

We now increase the parameters $\alpha_1$, $\alpha_2$, and $\alpha_3$ of the
previous experiment to consider a nonlinear Landau setting (see the beginning
of Section~\ref{sec:numexp} for the specific values of the parameters).
In this case, we set the space and velocity degrees of freedom to $64^3\times 128^3$,
while the tolerances for our algorithm are set to
$\epsilon_{\text{incr}}=10^{-4}$ and $\epsilon_{\text{decr}}=10^{-8}$.
Since we are considering a nonlinear regime, we expect the ranks needed by the
simulation to be higher compared to the linear case. For this reason,
the initial, minimum, and maximum ranks are set to
the values 25, 10, and 150, respectively. The final time is $T=40$ and the time
step size is chosen as $\tau=0.02$.
The obtained results are plotted in Figure~\ref{fig:flagship_nll}.
As anticipated, this example is more challenging due to its nonlinear nature,
and considerably higher ranks are required throughout the simulation.
With a total wall-clock time of approximately 43 minutes, our proposed
approach successfully reproduces the expected behavior of the electric energy:
after an initial damping phase, it increases at the known rate.

\subsubsection{Two-stream instability}\label{sec:ts}
\begin{figure}[htb]
  \centering
%
%
\begin{tikzpicture}

\begin{axis}[%
width=1.85in,
height=1in,
at={(1.62in,0.823in)},
scale only axis,
xmin=0,
xmax=40,
ymode=log,
ymin=0.0001,
ymax=10000,
yminorticks=true,
xlabel = {Time},
ylabel = {Electric energy},
axis background/.style={fill=white},
title={Electric energy evolution}
]
\addplot [color=red, forget plot]
  table[row sep=crcr]{%
0	0.581367687702657\\
0.0640000000000001	0.57862551963359\\
0.128	0.570454281344563\\
0.192	0.557018498266426\\
0.256	0.538588171425692\\
0.32	0.515532574804044\\
0.384	0.488311669828203\\
0.448	0.457465846162917\\
0.512	0.423603611458435\\
0.576000000000001	0.387387258703995\\
0.640000000000001	0.34951833484437\\
0.704000000000001	0.31072172434417\\
0.768000000000001	0.271729444638568\\
0.832000000000001	0.233264416556139\\
0.896000000000001	0.196024598586816\\
0.927999999999997	0.178071686719028\\
0.960000000000001	0.160667862425216\\
0.991999999999997	0.143886793939214\\
1.024	0.127797966338014\\
1.056	0.112466325865081\\
1.088	0.097951953173103\\
1.12	0.0843097673091042\\
1.152	0.0715892620660312\\
1.184	0.0598342761165058\\
1.216	0.0490827981280097\\
1.248	0.0393668078356898\\
1.28	0.0307121538206916\\
1.312	0.02313846850996\\
1.344	0.0166591206793322\\
1.376	0.0112812057610346\\
1.408	0.0070056143928013\\
1.44	0.00382697112910489\\
1.472	0.00173386154430496\\
1.504	0.000708926193291012\\
1.536	0.000729031050739917\\
1.568	0.00176547100490921\\
1.6	0.0037842049164402\\
1.632	0.00674612056519856\\
1.664	0.0106073276325281\\
1.696	0.0153194767035028\\
1.728	0.0208301021248281\\
1.76	0.0270829864198349\\
1.792	0.034018543843332\\
1.824	0.0415742205566049\\
1.856	0.0496849088171243\\
1.888	0.0582833725090391\\
1.92	0.0673006812895919\\
1.95200000000001	0.0766666505934221\\
1.984	0.0863102847214285\\
2.016	0.0961602202434568\\
2.048	0.10614516696429\\
2.08	0.116194343740243\\
2.112	0.126237906488434\\
2.144	0.13620736580225\\
2.176	0.146035991674029\\
2.20800000000001	0.155659202928718\\
2.27200000000001	0.17404401252867\\
2.33600000000001	0.190901743981958\\
2.40000000000001	0.205828307071777\\
2.46400000000001	0.21848413273704\\
2.528	0.22860094573112\\
2.592	0.235986617336796\\
2.656	0.240528017363507\\
2.72	0.24219183660569\\
2.784	0.241023402274075\\
2.848	0.237143558458936\\
2.912	0.23074373045862\\
2.976	0.222079335977747\\
3.04	0.211461748022778\\
3.104	0.199249048405634\\
3.168	0.185835831326475\\
3.232	0.171642330037212\\
3.296	0.157103152209391\\
3.36000000000001	0.142655915846743\\
3.456	0.122091996208046\\
3.584	0.0988176627499566\\
3.648	0.0897385297738332\\
3.68	0.0859659539837441\\
3.712	0.0827445387822202\\
3.744	0.0800996878633989\\
3.776	0.0780530730135379\\
3.808	0.0766225336241351\\
3.84	0.0758220009807193\\
3.872	0.075661447445078\\
3.904	0.0761468604879391\\
3.936	0.0772802413711587\\
3.968	0.0790596281243119\\
4	0.0814791423111583\\
4.032	0.0845290589378651\\
4.064	0.0881958987178358\\
4.128	0.0973083617734507\\
4.192	0.108638433040649\\
4.416	0.162070533472174\\
4.48	0.179970511999286\\
4.544	0.198386381613399\\
4.608	0.216940185269462\\
4.672	0.235253845294469\\
4.736	0.252957376487969\\
4.8	0.269696648075342\\
4.864	0.285140529658622\\
4.928	0.298987279382892\\
4.992	0.310970058232502\\
5.056	0.320861481856485\\
5.12	0.328477149780917\\
5.184	0.333678120483753\\
5.248	0.336372328780487\\
5.312	0.336514968518424\\
5.376	0.334107887973663\\
5.44	0.329198066957735\\
5.504	0.321875263085671\\
5.568	0.312268930013766\\
5.63200000000001	0.300544523041278\\
5.69600000000001	0.286899316466744\\
5.76000000000001	0.271557860016317\\
5.82400000000001	0.254767198133355\\
5.88800000000001	0.236791969682333\\
5.95200000000001	0.217909498991163\\
6.01600000000001	0.198404980862555\\
6.08000000000001	0.17856685113013\\
6.14400000000001	0.158682420623017\\
6.20800000000001	0.139033834908762\\
6.27200000000001	0.11989440611776\\
6.33600000000001	0.101525347367801\\
6.40000000000001	0.0841729252962073\\
6.46400000000001	0.0680660323718824\\
6.496	0.0605459326651891\\
6.52800000000001	0.0534141683891463\\
6.56	0.0466936295377526\\
6.59200000000001	0.0404058101833585\\
6.624	0.0345707766432969\\
6.65600000000001	0.0292071404210886\\
6.688	0.0243320356143438\\
6.72000000000001	0.0199611004776232\\
6.752	0.0161084628291893\\
6.78400000000001	0.0127867289956724\\
6.88	0.0061100476164994\\
6.91200000000001	0.00500734924576117\\
6.944	0.00447558851549042\\
6.97600000000001	0.00451817333284484\\
7.008	0.00513698895504915\\
7.04000000000001	0.00633240607793898\\
7.16800000000001	0.0168350789034512\\
7.20000000000001	0.0208668170384713\\
7.23200000000001	0.0254462227814965\\
7.264	0.0305634201501364\\
7.29600000000001	0.0362071384091972\\
7.328	0.0423647343804336\\
7.36000000000001	0.0490222167359212\\
7.392	0.0561642724440007\\
7.456	0.0718344185083106\\
7.52	0.0892273930455092\\
7.584	0.108176393355221\\
7.648	0.12849693644679\\
7.712	0.149988430105058\\
7.776	0.172435994955747\\
7.84	0.195612555223077\\
7.904	0.219281208749965\\
7.968	0.243197877920162\\
8.032	0.267114232624921\\
8.096	0.290780864674611\\
8.16	0.313950680617278\\
8.224	0.336382467672279\\
8.288	0.357844576205273\\
8.352	0.378118652131849\\
8.416	0.397003343765642\\
8.48	0.414317899105287\\
8.57600000000001	0.437010152758872\\
8.672	0.455411351794032\\
8.76800000000001	0.469265843423847\\
8.864	0.478487750396973\\
8.96000000000001	0.483167222167841\\
9.056	0.483568429567737\\
9.15200000000001	0.480118747698394\\
9.248	0.473389080538515\\
9.376	0.460512702748882\\
9.53600000000001	0.440744570586141\\
9.85600000000001	0.402619645779373\\
9.98400000000001	0.392013189078693\\
10.112	0.385300619145342\\
10.24	0.382494797102882\\
10.368	0.382994585908419\\
10.592	0.388167329552253\\
10.784	0.391357417441893\\
10.912	0.389922486885076\\
11.008	0.386058602636489\\
11.104	0.379472246209783\\
11.2	0.370132831497843\\
11.296	0.358297713839101\\
11.424	0.339659318002725\\
11.712	0.298150249942648\\
11.776	0.291980455306952\\
11.84	0.288030154562233\\
11.904	0.286846926938407\\
11.968	0.288978067805422\\
12.032	0.29495811713818\\
12.096	0.305296135776056\\
12.16	0.320463032957207\\
12.224	0.34087925498288\\
12.288	0.3669031433757\\
12.352	0.398820261446613\\
12.448	0.458166138084661\\
12.576	0.559045906835105\\
12.768	0.754985504862994\\
12.864	0.870651758524954\\
12.96	0.995784126137228\\
13.056	1.12810713843594\\
13.152	1.26501793726701\\
13.248	1.40368565278336\\
13.344	1.54116357997152\\
13.44	1.67450876642665\\
13.536	1.80090318755613\\
13.632	1.91777056052337\\
13.728	2.02288303430009\\
13.824	2.11445252703164\\
13.92	2.19120242769029\\
14.016	2.25241650866607\\
14.112	2.29796296170354\\
14.208	2.32829276021932\\
14.336	2.34687697539684\\
14.464	2.34401428537747\\
14.624	2.3180563671761\\
14.944	2.2373641594622\\
15.104	2.20770266500288\\
15.232	2.19959210407938\\
15.36	2.21060956385972\\
15.488	2.2446313038731\\
15.584	2.28706619206426\\
15.68	2.34508884132357\\
15.808	2.44797047895398\\
15.936	2.58124695718578\\
16.064	2.7460388693274\\
16.192	2.94345082595596\\
16.352	3.23806894744399\\
16.544	3.66610366738069\\
16.768	4.27712753497424\\
17.12	5.50785996693387\\
17.536	7.41803669741558\\
17.76	8.63617321017318\\
17.952	9.76226115668001\\
18.144	10.9393854526432\\
18.336	12.1408610572265\\
18.528	13.3404827497925\\
18.72	14.5165676079674\\
18.912	15.6558533465038\\
19.136	16.9370945260665\\
19.392	18.3664938419283\\
20.064	22.6002631803316\\
20.288	24.4871893775397\\
20.48	26.4373668148208\\
20.672	28.7702852836046\\
20.864	31.5527131892228\\
21.088	35.4448932988723\\
21.344	40.834289475871\\
21.792	52.8617566888411\\
22.208	67.0039776263701\\
22.528	79.706413605246\\
22.816	92.386857072244\\
23.104	106.194807302863\\
23.424	122.90355412829\\
23.808	145.195759779426\\
24.96	237.998023737422\\
25.344	284.095032232757\\
25.76	347.334049461014\\
26.4	478.173953753172\\
27.008	646.027928105087\\
27.456	799.457891095796\\
27.872	966.228379647629\\
28.288	1158.37004905512\\
28.736	1396.723672104\\
29.216	1692.40904225246\\
29.664	2008.47628920413\\
30.048	2308.53065461819\\
30.4	2602.24802573832\\
30.72	2878.21208287493\\
31.008	3127.17507574846\\
31.296	3369.71143870736\\
31.584	3598.8125332798\\
31.84	3785.7088994857\\
32.096	3952.36412283327\\
32.352	4095.38338739627\\
32.64	4225.32400615801\\
32.9279999999999	4321.65136518959\\
33.2159999999999	4385.23002358944\\
33.5359999999999	4420.58611952059\\
33.8879999999998	4422.61143034395\\
34.3039999999998	4386.75363919554\\
34.8159999999997	4304.44976871359\\
35.4559999999996	4167.56674767024\\
36.2879999999995	3961.38132458343\\
37.1839999999995	3718.46561136488\\
38.3679999999993	3418.68984659609\\
38.8799999999993	3325.9734166854\\
39.3599999999992	3269.12004495759\\
39.8399999999992	3241.18089539617\\
39.9679999999991	3237.92106050292\\
};
\end{axis}

\end{tikzpicture}%
  \input{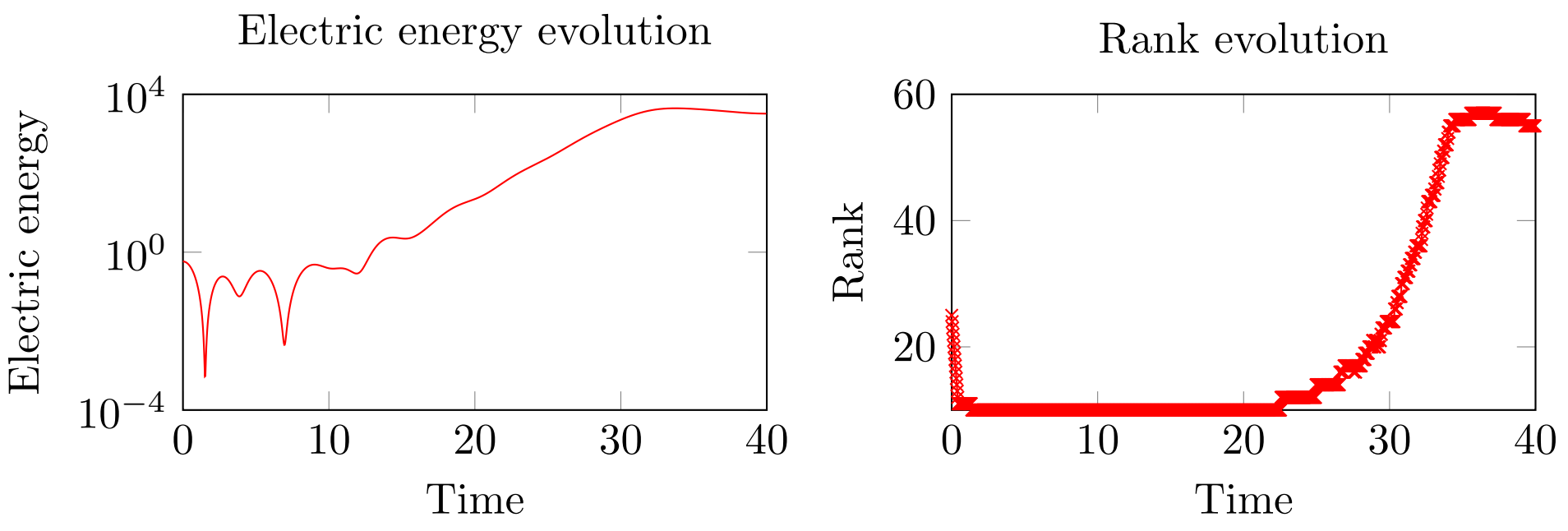}
  \caption{Results for the $128^3\times 128^3$ two-stream instability benchmark example
  (see Section~\ref{sec:ts} for more details).
  Left: electric energy evolution.
  Right: rank evolution.}
  \label{fig:flagship_ts}
\end{figure}

We now consider a different benchmark example: a two-stream
instability simulation (see equation~\eqref{eq:ts} for the precise initial
condition). The degrees of freedom are $128^3$ for both the space and velocity
variables, while the tolerances for the algorithm are
$\epsilon_{\text{incr}}=10^{-5}$ and $\epsilon_{\text{decr}}=10^{-8}$.
The time step size of the numerical scheme is $\tau=0.032$ with the final time
set to $T=40$.
The initial, minimum, and maximum ranks are 25, 10, and 100, respectively.
The results of the simulation are presented in Figure~\ref{fig:flagship_ts}.
Here, too, the results are as expected: the rank remains essentially constant
at a low value of 10 during the linear phase and then progressively increases
as the nonlinear phase begins.
The wall-clock time of the overall simulation is approximately 19 minutes.
Note that, by employing our adaptive DLR approach, we are able
to efficiently perform the simulation using a collection of objects with a
reserved memory size of $128^3\times 100$. In contrast, a full-rank approach would
require storing objects of size $128^6$, making the computation infeasible due
to the associated memory requirements.

\subsubsection{Bump-on-tail problem}\label{sec:bot}
\begin{figure}[htb]
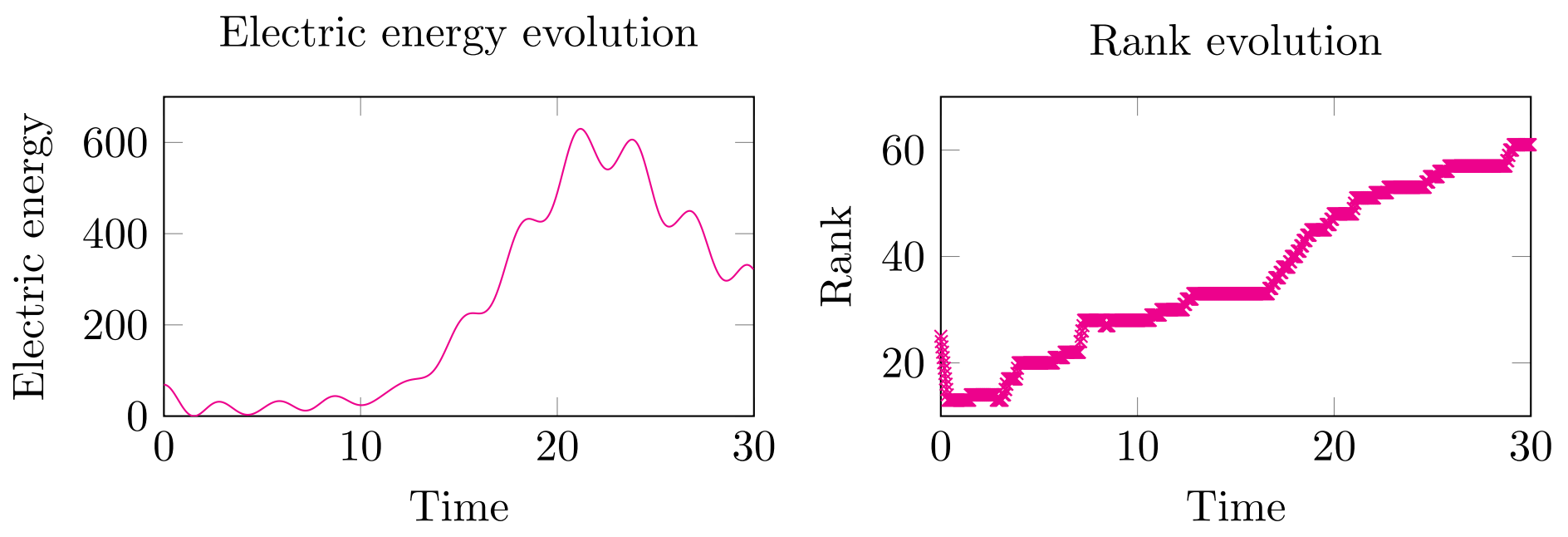

  \centering
  \input{plot_ee_bot.tex}
  \input{plot_rank_bot.tex}
  \caption{Results for the $128^3\times 128^3$ bump-on-tail benchmark example
  (see Section~\ref{sec:bot} for more details).
  Left: electric energy evolution.
  Right: rank evolution.}
  \label{fig:flagship_bot}
\end{figure}

We conclude the tested benchmark examples with the bump-on-tail problem defined
by the initial condition~\eqref{eq:bot}.
In this case, we set the space and velocity degrees of freedom to $128^3$ in both
space and velocity, while the tolerances of the adaptive PS algorithm are set to
$\epsilon_{\text{incr}}=10^{-4}$ and $\epsilon_{\text{decr}}=10^{-8}$.
The initial rank is set to 25, the minimum rank to 10, and the maximum rank to
100. We select the final time $T=30$ and the time step size $\tau=0.03$.
The results of the simulation are depicted in Figure~\ref{fig:flagship_bot} and
again agree with the expectations. In particular, the
evolution of the electric energy is accurately captured, closely matching results
previously reported in the literature.
The wall-clock time of the overall simulation is approximately 26 minutes.

\section{Conclusions}\label{sec:concl}

In this work, we presented an adaptive dynamical low-rank approach for the
efficient solution of the six-dimensional Vlasov--Poisson system.
Built on projector-splitting-type integrators, our strategy exploits
physically relevant quantities, such as the electric field and energy, to
dynamically adjust the rank of the approximation throughout the simulation.
To this end, we extended the dynamical low-rank framework Ensign,
in which the proposed methodology is implemented.
The key advantage of the proposed approach is its suitability for high-performance
computing architectures, particularly GPUs. In fact, compared, e.g., with
state-of-the-art adaptive Basis-Update \& Galerkin
integrators, our method operates with smaller intermediate low-rank objects,
resulting in a reduced memory footprint. This is particularly relevant in GPU-based
simulations, where memory capacity can be a limiting factor, and allows our
approach to outperform established techniques in the considered HPC setting.
We extensively tested our approach on a broad range of standard 6D benchmarks,
including linear and nonlinear Landau damping, the two-stream instability,
and the bump-on-tail instability. The results demonstrate that the method
robustly captures the relevant plasma dynamics while maintaining a compact
low-rank representation.

\section*{Acknowledgments}
Fabio Cassini is a member of the Gruppo Nazionale
Calcolo Scientifico-Istituto Nazionale di Alta Matematica (GNCS-INdAM).

This research was funded [in part] by the Austrian Science Fund (FWF), \url{https://doi.org/10.55776/PAT2937525}.

\bibliographystyle{plain}
\bibliography{CE26}

@book{B13,
  author={Bittencourt, J. A.},
  title={Fundamentals of plasma physics},
  year={2013},
  publisher={Springer New York, NY},
  edition = {third},
  doi = {https://doi.org/10.1007/978-1-4757-4030-1}
}

@article{DP14,
  author = {Dimarco, G. and Pareschi, L.},
  title = {Numerical methods for kinetic equations},
  journal = {Acta Numer.},
  year = {2014},
  pages = {369--520},
  volume = {23},
  doi = {https://doi.org/10.1017/S0962492914000063}
}

@article{KL07,
  author = {Koch, O. and Lubich, C.},
  title = {Dynamical {L}ow-{R}ank approximation},
  journal = {SIAM J. Matrix. Anal. Appl.},
  volume = {29},
  pages = {434--454},
  year = {2007},
  doi = {https://doi.org/10.1137/050639703}
}

@article{CE22,
  author = {Cassini, F. and Einkemmer, L.},
  title = {Efficient 6{D} {V}lasov simulation using the dynamical low-rank
  framework \texttt{Ensign}},
  journal = {Comput. Phys. Commun.},
  volume = {280},
  pages = {108489},
  year = {2022},
  doi = {https://doi.org/10.1016/j.cpc.2022.108489}
}

@article{EL18,
  author = {Einkemmer, L. and Lubich, C.},
  title = {A {L}ow-{R}ank {P}rojector-{S}plitting {I}ntegrator for the {V}lasov--{P}oisson {E}quation},
  journal = {SIAM J. Sci. Comput.},
  volume = {40},
  number = {5},
  pages = {B1330--B1360},
  year = {2018},
  doi = {https://doi.org/10.1137/18M116383X}
}

@article{HNS23,
  author = {Hochbruck, M. and Neher, M. and Schrammer, S.},
  title = {Rank-adaptive dynamical low-rank integrators for first-order and second-order matrix differential equations},
  journal = {BIT Numer. Math.},
  volume = {63},
  number = {9},
  pages = {1--24},
  year = {2023},
  doi = {https://doi.org/10.1007/s10543-023-00942-6}
}

@article{CEKL24,
  author = {Ceruti, G. and Einkemmer, L. and Kusch, J. and Lubich, C.},
  title = {A robust second-order low-rank BUG integrator based on the midpoint rule},
  journal = {BIT Numer. Math.},
  volume = {64},
  number = {30},
  pages = {1--19},
  year = {2024},
  doi = {https://doi.org/10.1007/s10543-024-01032-x}
}

@article{LO14c,
  author = {Lubich, C. and Oseledets, I. V.},
  title = {A projector-splitting integrator for dynamical low-rank approximation},
  journal = {BIT Numer. Math.},
  volume = {54},
  pages = {171--188},
  year = {2014},
  doi = {https://doi.org/10.1007/s10543-013-0454-0}
}

@article{EOP20,
  author = {Einkemmer, L. and Ostermann, A. and Piazzola, C.},
  title = {A low-rank projector-splitting integrator for the {V}lasov--{M}axwell equations with divergence correction},
  journal = {J. Comput. Phys.},
  volume = {403},
  pages = {109063},
  year = {2020},
  doi = {https://doi.org/10.1016/j.jcp.2019.109063}
}

@article{PMF20,
  author = {Peng, Z. and McClarren, R. G. and Frank, M.},
  title = {A low-rank method for two-dimensional time-dependent radiation transport calculations},
  journal = {J. Comput. Phys.},
  volume = {421},
  pages = {109735},
  year = {2020},
  doi = {https://doi.org/10.1016/j.jcp.2020.109735}
}

@article{PEL23,
  author = {Prugger, M. and Einkemmer, L. and Lopez, C. F.},
  title = {A dynamical low-rank approach to solve the chemical master equation for biological reaction networks},
  journal = {J. Comput. Phys.},
  volume = {489},
  pages = {112250},
  year = {2023},
  doi = {https://doi.org/10.1016/j.jcp.2023.112250}
}

@article{EMP24,
  author = {Einkemmer, L. and Mangott, J. and Prugger, M.},
  title = {A low-rank complexity reduction algorithm for the high-dimensional kinetic chemical master equation},
  journal = {J. Comput. Phys.},
  volume = {503},
  pages = {112827},
  year = {2024},
  doi = {https://doi.org/10.1016/j.jcp.2024.112827}
}

@article{CL22,
  author = {Ceruti, G. and Lubich, C.},
  title = {An unconventional robust integrator for dynamical low-rank approximation},
  journal = {BIT Numer. Math.},
  volume = {62},
  pages = {23--44},
  year = {2022},
  doi = {https://doi.org/10.1007/s10543-021-00873-0}
}

@article{KCEF22,
  author = {Kusch, J. and Ceruti, G. and Einkemmer, L. and Frank, M.},
  title = {Dynamical low-rank approximation for Burgers' equation with uncertainty},
  journal = {Int. J. Uncertain. Quantif.},
  volume = {12},
  issue = {5},
  pages = {1--21},
  year = {2022},
  doi = {https://doi.org/10.1615/Int.J.UncertaintyQuantification.2022039345}
}

@article{EKKMQ25,
  author = {Einkemmer, L. and Kormann, K. and Kusch, J. and McClarren, R. G.
            and Qiu, J.-M.},
  title = {A review of low-rank methods for time-dependent kinetic equations},
  journal = {J. Comput. Phys.},
  volume = {538},
  pages = {114191},
  year = {2025},
  doi = {https://doi.org/10.1016/j.jcp.2025.114191}
}

@article{CLS23,
  author = {Ceruti, G. and Lubich, C. and Sulz, D.},
  title = {Rank-adaptive time integration of tree tensor networks},
  journal = {SIAM J. Numer. Anal.},
  volume = {61},
  issue = {1},
  pages = {194-222},
  year = {2023},
  doi = {https://doi.org/10.1137/22M1473790}
}

@article{CKL24,
  author = {Ceruti, G. and Kusch, J. and Lubich, C.},
  title = {A parallel rank-adaptive integrator for dynamical low-rank approximation},
  journal = {SIAM J. Sci. Comput.},
  volume = {43},
  issue = {3},
  pages = {B205-B228},
  year = {2024},
  doi = {https://doi.org/10.1137/23M1565103}
}

@article{HW22,
  author = {Hu, J. and Wang, Y.},
  title = {An adaptive dynamical low rank method for the nonlinear Boltzmann equation},
  journal = {J. Sci. Comput.},
  volume = {92},
  number = {75},
  pages = {1--24},
  year = {2022},
  doi = {https://doi.org/10.1007/s10915-022-01934-4}
}

@misc{cublas,
  title={{cuBLAS documentation}},
  author = {NVIDIA Corporation},
  howpublished = {\url{https://docs.nvidia.com/cuda/cublas/index.html}},
  year = {2026}
}

@inproceedings{ZWZ12,
  title={Model-driven {L}evel 3 {BLAS} {P}erformance {O}ptimization on {L}oongson 3{A} {P}rocessor},
  author={Zhang, X. and Wang, Q. and Zhang, Y.},
  booktitle={2012 IEEE 18th international conference on parallel and distributed systems},
  pages={684--691},
  year={2012},
  publisher={IEEE},
  doi = {https://doi.org/10.1109/ICPADS.2012.97}
}

@inproceedings{QZZY13,
  title={{AUGEM}: automatically generate high
  performance dense linear algebra kernels
  on x86 {CPU}s},
  author={Qian, W. and Zhang, X. and Zhang, Y. and Yi, Q.},
  booktitle={2013 International Conference on High Percormance Computing, Networking, Storage and Analysis},
  pages={25},
  year={2013},
  publisher={ACM},
  doi = {https://doi.org/10.1145/2503210.2503219}
}

@article{FJ05,
  author = {M. Frigo and Johnson, S. G.},
  title = {The Design and Implementation of {FFTW3}},
  journal = {Proceedings of the IEEE},
  year = {2005},
  volume = {93},
  number = {2},
  pages = {216--231},
  doi = {https://doi.org/10.1109/JPROC.2004.840301}
}

@article{ceruti2022rank,
  title={A rank-adaptive robust integrator for dynamical low-rank approximation},
  author={Ceruti, G. and Kusch, J. and Lubich, C.},
  journal={BIT Numer. Math.},
  volume={62},
  number={4},
  pages={1149--1174},
  year={2022},
  doi = {https://doi.org/10.1007/s10543-021-00907-7}
}

@article{ehrlacher2017,
  title={A dynamical adaptive tensor method for the {V}lasov--{P}oisson system},
  author={Ehrlacher, V. and Lombardi, D.},
  journal={J. Comput. Phys.},
  volume={339},
  pages={285--306},
  year={2017},
  doi = {https://doi.org/10.1016/j.jcp.2017.03.015}
}

@article{guo2022,
  title={A low rank tensor representation of linear transport and nonlinear {V}lasov solutions and their associated flow maps},
  author={Guo, W. and Qiu, J.-M.},
  journal={J. Comput. Phys.},
  volume={458},
  pages={111089},
  year={2022},
  doi = {https://doi.org/10.1016/j.jcp.2022.111089}
}

@article{einkemmer2023semi,
  title={Semi-{L}agrangian 4d, 5d, and 6d kinetic plasma simulation on large-scale {GPU}-equipped supercomputers},
  author={Einkemmer, L. and Moriggl, A.},
  journal={Int. J. High Perform. Comput. Appl.},
  volume={37},
  number={2},
  pages={180--196},
  year={2023},
  doi = {https://doi.org/10.1177/10943420221137599}
}

@article{kormann2019,
  title={A massively parallel semi-{L}agrangian solver for the six-dimensional {V}lasov--{P}oisson equation},
  author={Kormann, K. and Reuter, K. and Rampp, M.},
  journal={Int. J. High Perform. Comput. Appl.},
  volume={33},
  number={5},
  pages={924--947},
  year={2019},
  doi = {https://doi.org/10.1177/1094342019834644}
}

@article{burau2010,
  title={{PIC}on{GPU}: {A} fully relativistic particle-in-cell code for a {GPU} cluster},
  author={Burau, H. and Widera, R. and H{\"o}nig, W. and Juckeland, G. and Debus, A. and Kluge, T. and Schramm, U. and Cowan, T.E. and Sauerbrey, R. and Bussmann, M.},
  journal={IEEE Trans. Plasma Sci.},
  volume={38},
  number={10},
  pages={2831--2839},
  year={2010},
  doi = {https://doi.org/10.1109/TPS.2010.2064310}
}

@article{hauck2023,
  title={A {P}redictor-{C}orrector {S}trategy for {A}daptivity in {D}ynamical {L}ow-{R}ank {A}pproximations},
  author={Hauck, C. D. and Schnake, S.},
  journal={SIAM J. Matrix Anal. Appl.},
  volume={44},
  number={3},
  pages={971--1005},
  year={2023},
  doi = {https://doi.org/10.1137/22M1519493}
}

\end{document}